\documentclass[11pt]{article}  

\usepackage{epsfig}
\usepackage{amsfonts, amsmath, amssymb, color, mathrsfs, mathabx}  
 \usepackage{amsthm}
 \usepackage{placeins}

\newtheorem{thm}{Theorem}
\newtheorem{prop}[thm]{Proposition} 

\newcommand{\mc}[1]{\mathcal{#1}}
\newcommand{\mb}[1]{\mathbb{#1}}

\newcommand{\arrow}{\xrightarrow{}}

\usepackage[subpreambles=true]{standalone}
\usepackage{import}
\usepackage{subfiles} 

\begin{document}
\title{Bifurcation analysis of a conceptual 
  model for the Atlantic Meridional
	Overturning Circulation}
\author{John Bailie and Bernd Krauskopf}
\maketitle

\begin{abstract}
\noindent
The Atlantic Meridional Overturning Circulation (AMOC) distributes heat and salt into the Northern Hemisphere via a warm surface current toward the subpolar North Atlantic, where water sinks and returns southwards as a deep cold current. There is substantial evidence that the AMOC has slowed down over the last century. We introduce a conceptual box model for the evolution of salinity and temperature on the surface of the North Atlantic Ocean, subject to the influx of meltwater from the Greenland ice sheets. Our model, which extends a model due to Welander, describes the interaction between a surface box and a deep-water box of constant temperature and salinity, which may be convective or non-convective, depending on the density difference. Its two main parameters $\mu$ and $\eta$ describe the influx of freshwater and the threshold density between the two boxes, respectively.

We use tools from bifurcation theory to analyse two cases of the model: the limiting case of instantaneous switching between convective or non-convective interaction, where the system is piecewise-smooth (PWS), and the full smooth model with more gradual switching. For the PWS model we perform a complete bifurcation analysis by deriving analytical expressions for all bifurcations. The resulting bifurcation diagram in the $(\mu,\eta)$-plane identifies all regions of possible dynamics, which we show as phase portraits --- both at typical parameter points, as well as at the different transitions between them. We also present the bifurcation diagram for the case of smooth switching and show how it arises from that of the PWS case. In this way, we determine exactly where one finds bistability and self-sustained oscillations of the AMOC in both versions of the model. In particular, our results show that oscillations between temperature and salinity on the surface North Atlantic Ocean disappear completely when the transition between the convective and non-convective regimes is too slow. 
\end{abstract}

\section{Introduction}
The Atlantic Meridional Overturning Circulation (AMOC) is a large conveyor belt of water that spans the entire Atlantic Ocean. Light surface currents transport relatively warm and saline waters northward to high latitudes. Here, the water becomes denser, leading to downward convection and mixing with the deep ocean, and subsequent formation of deepwater masses.  A deep current then transports this water back to lower latitudes, where it upwells to the surface, thus closing the circulation loop \cite{kuhlbrodt2007driving}. The strength of the AMOC is governed by the interplay between two proposed upwelling mechanisms \cite{speer2000diabatic, sloyan2001southern}. The first perspective is that turbulent mixing across surfaces of equal density results in the upwelling of deepwater to the surface ocean in low latitudes \cite{jeffreys1925fluid, sandstrom1916meteorologische}. The second perspective suggests that strong circumpolar winds induce upwelling in the South Atlantic Ocean \cite{toggweiler1998ocean}.  Regardless of the mechanism, the process of deepwater formation is crucial in determining the shape and strength of the associated return current --- making it a critical factor for the stability of the AMOC.

\begin{figure}[t!]
	\centering
	\includegraphics{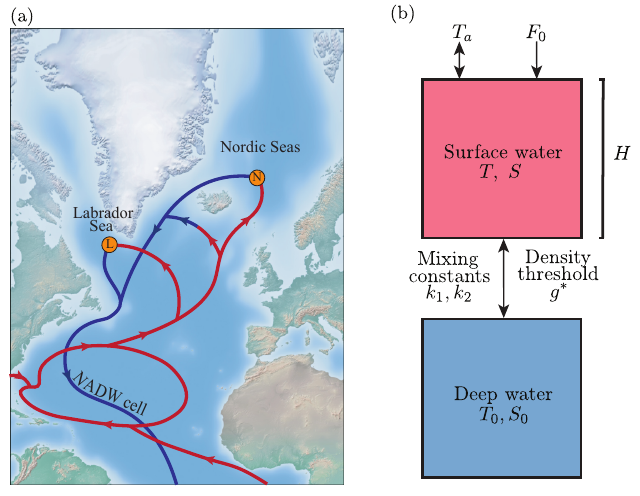}
	\caption{Panel~(a) shows a simplified sketch of the North Atlantic component of the AMOC, inspired by \cite{srokosz2012past}. The surface flow is displayed in red, and the NADW overturning cell in blue. Deepwater formation sites are denoted $\mathrm{L}$ and $\mathrm{N}$ in the Labrador sea and the Nordic seas, respectively. Panel~(b) shows the two-box model setup for the interaction between surface water and cold deep water at the sites $\mathrm{L}$ and $\mathrm{N}$.}
	\label{fig:intro_fig1}
\end{figure}

This paper focuses on the deepwater formation sites in the North Atlantic. Specifically, as illustrated in Figure~\ref{fig:intro_fig1}(a), the convection of highly saline water from the surface to the deep ocean in the Labrador and Nordic seas forms the North Atlantic Deep Water (NADW); it has an associated return current referred to as the NADW overturning cell. Several climate processes, such as salt rejection and atmospheric cooling, facilitate this convection by preconditioning the subpolar North Atlantic to have relatively high salinity \cite{marsh2007stability}. Furthermore, an advective process transports saline water to the North Atlantic and, thus, stimulates the convection and formation of the NADW \cite{weijer2019stability}. An inherent negative feedback loop forms: weaker convection results in a smaller NADW and, consequently, a weaker overturning cell. This weaker cell then advects less salt to the North Atlantic, further weakening the convection.

Evidence from proxy and sea surface temperature measurements indicate that the AMOC has weakened over the twentieth century \cite{rahmstorf2015exceptional}. There was also a particularly abrupt change in overturning strength during the 1970s \cite{dima2010evidence} attributed to a large-scale influx of fresh water into the North Atlantic; this is known as the Great Salinity Anomaly and is linked to Arctic sea-ice export \cite{dickson1988great}. A weakened AMOC has significant consequences on the Earth's climate system since it leads to reduced northern heat transport, which lowers the oceanic and atmospheric temperature in the Northern hemisphere \cite{jackson2015global} via a weakening or even shutdown of the northern deepwater formation. Some significant implications drawn from simulations are a widespread cooling in Europe \cite{rahmstorf2015exceptional}, the possible collapse of the North Atlantic plankton stocks \cite{schmittner2005decline}, and a rise in the sea level \cite{gregory2005model}. As a result of external environmental factors, the AMOC is likely to weaken further, and a complete shutdown of the deep water formation in the Labrador Sea is a possibility \cite{rahmstorf2015exceptional}. In particular, meltwater from the melting Greenland ice sheets contributes to a large influx of freshwater into the subpolar North Atlantic \cite{nghiem2012extreme}. As freshwater is strictly non-saline, it dilutes the ocean surface water by lowering its salinity, thus, inhibiting the deep-water formation and, hence, the NADW overturning cell strength. 

Of particular interest in this context is the modelling of the underlying deep water formation itself --- with the aim of understanding the possible long-term behaviour of the AMOC in response to freshwater influx. Climate models form a hierarchy of complexity, and the choice of model depends on the nature of the question that is being asked. We study here a conceptual model of low complexity for investigating the AMOC in regard to deep water formation --- specifically, from the class of box models that consider only a few variables in a relatively small number of interacting boxes, each representing a body of water of concern. While they are not designed to be used for prediction, box models are simple enough to be amenable to mathematical analysis, including with tools from dynamical systems theory \cite{dijkstra2005nonlinear}. 

The stability of the AMOC was first investigated by Stommel with a two-box model \cite{stommel1961thermohaline}; it considers the circulation between a subtropical box and a subpolar box, where a capillary flow represents the advection of water between the two boxes. Stommel's model features three qualitatively different regimes. In the first regime, the AMOC is driven by salinity differences between the boxes, and surface currents move water toward the equator. Temperature differences are the main driver in the second regime, and the surface currents move water toward the poles. The final regime features bistability, where the AMOC may tip to either of the described equilibrium states. Stommel laid the foundation for several advective models, which add more boxes and physical processes; see, for example, \cite{rooth1982hydrology, welander1986thermohaline,nkk_3box}. 

A two-box model presented by Welander \cite{Walender1982} attempts to describe self-sustained oscillations of temperature and salinity on the ocean surface in the presence of external forcing. The boxes interact by exchanging heat and salt via a mixing process. When the water in the surface box is sufficiently dense, the mixing is convective (strong). When the water in the boxes has comparable density, on the other hand, the mixing is non-convective (weak) and may happen via several climate processes, such as double-diffusion \cite{huppert1981double}. In this setup, an atmospheric basin with fixed properties interacts with the surface box, which is modelled by Newton's transfer law. The model by Welander is described for two cases: when the transition between convective and non-convective mixing is modelled as a continuous change and, alternatively, when it is instantaneous and discontinuous. In both cases, self-sustained oscillations are observed, which are characterised by a convective and a non-convective phase. Welander's model was re-examined by Leifeld \cite{leifeld2016nonsmooth} with the aim of formalising the previous analysis by using a modern approach of piecewise-smooth (PWS) dynamical systems. They undertook a preliminary stability analysis and made a first comparison between the smooth and non-smooth models; however, this work falls short of describing the full bifurcation picture and, to the best of our knowledge, there is as yet no complete analysis of the Welander model, nor any closely related models. 

\subsection{The adjusted Welander model}

We take this as the starting point of our study of an \emph{adjusted Welander model} that also considers the impact of a freshwater influx into the North Atlantic ocean. Following on from work in \cite{PCessi1996}, where the external forcing enters in the form of Newton's transfer law, we consider here a direct freshwater flux that dilutes the salinity in a surface ocean box at the North Atlantic, which is coupled to a box of deep water of constant lower temperature and salinity. As is illustrated by the schematic in Figure~\ref{fig:intro_fig1}(b), the model takes the form of a planar system of ordinary differential equations for temperature $T$ and salinity $S$ in the surface ocean box, which is given by 
\begin{equation}
\begin{aligned}
	&\frac{dT}{dt} = -\gamma(T - T_a) - k_\varepsilon(\rho)(T - T_0), \label{eq:adjusted_welander}
	\\
	&\frac{dS}{dt} = \frac{F_0}{H}S_0-k_\varepsilon(\rho)(S - S_0).
\end{aligned}
\end{equation}
The atmosphere externally drives the surface ocean box to a thermal equilibrium $T_a$ at rate $\gamma$, which is Newton's transfer law. The salinity, on the other hand, is directly forced by the freshwater flux $F_0$ at the rate $\frac{F_0}{H}S_0$, where $H$ is the depth of the surface ocean box. Moreover, $T_0$ and $S_0$ are the (fixed) temperature and salinity of the deep-ocean box that drive $T$ and $S$, respectively, as given by the convective exchange function $k_\varepsilon$. This function determines the coefficient for Newton's transfer law and takes as its argument the density $\rho$ of the surface ocean box given (in linear approximation) by 
\begin{align}
	\label{eq:eq_of_state}
	\frac{\rho}{\rho_0} = 1 + \alpha_S(S - S_0) - \alpha_T(T - T_0).
\end{align}
Here, the constant $\rho_0$ is the density of the bottom box, and the coefficients $\alpha_S$ and $\alpha_T$ are, respectively, the saline expansion and thermal compression constants \cite{dijkstra2005nonlinear}. 

\begin{figure}[t!]
	\centering
	\includegraphics{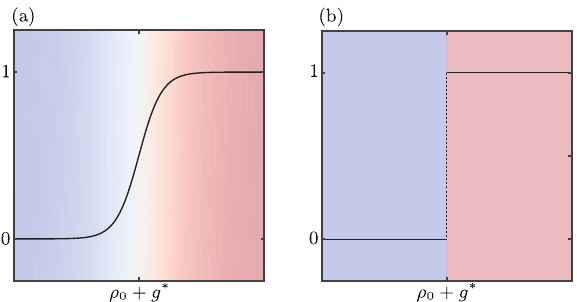}
	\caption{The switching functions $\mc{H}_{\varepsilon}(\rho-\rho_0-g^{*})$ with $\varepsilon = 0.1$ in panel~(a) and $\varepsilon =0$ in panel~(b). Blue shading indicates mixing is mainly non-convective, and red shading that mixing is mainly convective. In the smooth transition region in panel~(a), the colour changes from blue, via white, to red.}
	\label{fig:intro_fig2}
\end{figure}

The convective exchange function is a key ingredient in~\eqref{eq:adjusted_welander} and describes the transition between the two regimes when the vertical mixing between the two boxes is non-convective at rate $k_1 > 0$ and when it is convective at rate $k_2 > k_1$. This transition is modeled in its general form as 
\begin{align}
	\label{eq:convective_exchange}
	k_\varepsilon(\rho) = k_1 + \mc{H}_\varepsilon(\rho - \rho_0 - g^*)(k_2 - k_1),
\end{align}
where $\mc{H}_\varepsilon$ is a suitable switching function from zero to one, whose switching time depends on the switching-time parameter $\varepsilon$. Hence, when the density difference $\rho - \rho_0$ is (sufficiently) greater than the density threshold $g^*$, mixing between the boxes is mainly convective; on the other hand, it is mainly non-convective when $\rho - \rho_0$ is (sufficiently) smaller than $g^*$. Different switching functions have been used in the literature \cite{leifeld2016nonsmooth, welander1986thermohaline}, including those based on the $\arctan$ function. In this paper, we define $\mc{H}_\varepsilon$ as
\begin{align}
	\label{eq:transition_function}
	&\mc{H}_\varepsilon(u) = \frac{1}{2}\left( 1 + \tanh\left(\frac{u}{\varepsilon}\right)\right).
\end{align}
Note that $\mc{H}_\varepsilon$ has a switching time of order $\varepsilon>0$ and maximal rate of switching $\frac{1}{\varepsilon}$ given by the derivative of $\mc{H}_\varepsilon$ at zero. Moreover, the limiting case for $\varepsilon \arrow 0$ is an instantaneous switch, represented by the Heaviside function $\mc{H}_0$. Figure~\ref{fig:intro_fig2} shows the resulting convective exchange functions $k_\varepsilon(\rho)$ from~\eqref{eq:convective_exchange} for $\varepsilon=0.1$ and $\varepsilon=0$.

The adjusted Welander model in the form~\eqref{eq:adjusted_welander} has ten parameters, making a direct analysis impractical. The first step in our analysis is to non-dimensionalise the system by introducing rescaled temperature, salinity and time 
\begin{equation}
x = \frac{T - T_0}{T_a - T_0}, \quad 
y = \frac{\alpha_S(S - S_0)}{\alpha_T(T_a - T_0)},	\quad 
\tau = \gamma t, 
\end{equation} 
and parameters 
\begin{equation}
\kappa_i = \frac{k_i}{\gamma}, \quad
\mu = \frac{F_0S_0\alpha_S}{\gamma\alpha_T(T_a - T_0)H} \quad 
\eta = \frac{g^{*}(\kappa_2 - \kappa_1)}{\gamma\alpha_T(T_a - T_0)\rho_0}.
\end{equation} 
This transforms~\eqref{eq:adjusted_welander} into 
\begin{equation}
\begin{aligned}
	&\dot{x} = 1 - (1 + \kappa_1 + \mc{H}_\varepsilon(y - x - \eta)(\kappa_2 - \kappa_1)) x,
	\\
	&\dot{y} = \mu - (\kappa_1 + \mc{H}_\varepsilon(y - x - \eta)(\kappa_2 - \kappa_1))y,
	\label{eq:adjusted_welander_non_dim}
\end{aligned}
\end{equation}
where the dot represents the derivative with respect to the rescaled time.

\subsection{Outline of the work}

The adjusted Welander model in the form~\eqref{eq:adjusted_welander_non_dim} is our central object of study. We first perform in Section~\ref{section:pws} a (non-smooth) bifurcation analysis for the limiting case of system~\eqref{eq:adjusted_welander_non_dim} with the Heaviside switching function $\mc{H}_0$. Specifically, we determine and catalogue all of the possible dynamics, by presenting analytical expressions for all codimension-one and codimension-two bifurcations; the corresponding proofs and derivations can be found in Appendix A. The rescaled freshwater flux $\mu$ and density threshold $\eta$ are the bifurcation parameters, and we show the complete bifurcation diagram in the $(\mu,\eta)$-plane for a reasonable choice of the vertical mixing coefficients $0<\kappa_1<\kappa_2$. Moreover, we present representative phase portraits in the $(x,y)$-plane for all open regions, of which there are eight, for the different types of transitions of codimension one between them, as well as at the five codimension-two points that organise the bifurcation diagram. In particular, we identify the parameter regime where the system exhibits bistability between states, where deep-water convection is either substantial or shut down, which is characteristic behaviour of several models of different complexity \cite{weijer2019stability}. Moreover, we determine the parameter regime with self-sustained relaxation-type oscillations that have been observed in \cite{PCessi1996} and also in Welander's original work \cite{welander1986thermohaline}. In addition to this earlier work, we determine this region analytically and clarify the nature of the different possible transitions to/from this oscillatory regime. Our results also show that the bifurcation diagram in the $(\mu,\eta)$-plane is toplogically the same for any fixed $0<\kappa_1<\kappa_2$. 

Section~\ref{section:smooth} is then concerned with the smooth case of system~\eqref{eq:adjusted_welander_non_dim} with $\mc{H}_\varepsilon$ for $\varepsilon > 0$. Here we first present the bifurcation diagram in the $(\mu,\eta)$-plane for $\varepsilon = 0.1$ (for the same choice of $0<\kappa_1<\kappa_2$). This requires computing the relevant bifurcation curves by making use of established bifurcation theory \cite{kuznetsov1998elements} in conjunction with the continuation software package AUTO-07p \cite{auto_cont}. We focus here on the main parameter regimes, especially those that feature bistability and self-sustained oscillations, for which we show representative phase portraits. We then present a partial bifurcation analysis in $(\mu,\eta,\varepsilon)$-space that clarifies the convergence of the bifurcation diagram in the $(\mu,\eta)$-plane as $\varepsilon$ approaches $0$. Moreover, we show that there is a codimension-three bifurcation at a quite low value of the switching-time parameter $\varepsilon$, at which the region with self-sustained oscillations completely disappears from the $(\mu,\eta)$-plane. In other words, the switching between the regimes with strong convective mixing and with weak non-convective mixing needs to be sufficiently fast for relaxation-type oscillations to occur in the adjusted Welander model~\eqref{eq:adjusted_welander_non_dim}. In particular, this shows the relevance of the non-smooth limiting system with $\mc{H}_0$ for explaining this oscillatory behaviour. 

In the final Section~\ref{section:disscussion} we summarise our findings, briefly discuss their significance for the dynamics of AMOC, and point out some direction for future work.

\section{Bifurcation analysis of the PWS model for $\varepsilon = 0$}
\label{section:pws}

In the limiting case of an instantaneous transition with the transition function $\mc{H}_0$ from~\eqref{eq:transition_function}, system~\eqref{eq:adjusted_welander_non_dim} reduces to the piecewise-smooth linear Filippov system
\begin{align}
	\label{model:PWS}
	\binom{\dot{x}}{\dot{y}} = 
	\left\{
	\begin{matrix}
		&f_1(x, y), \ \ \ y < x + \eta
		\\
		&f_2(x, y), \ \ \ y > x + \eta
	\end{matrix}\right.
\end{align}
with
\begin{align}
	\label{model:PWSfi}
f_i(x,y) = \binom{1 - (1+\kappa_i)x}{\mu - \kappa_iy}.
\end{align}
The switching manifold 
\begin{align}
	\Sigma = \{(x,y)\in\mb R^2, \ \ y = x + \eta\}
	\label{eq: switching_manifold}
\end{align} 
is a straight line that partitions the phase space of~\eqref{model:PWS} into the open regions
\begin{align}
	\label{eq:R1}
	&R_1 = \{(x,y)\in\mb R^2, \  \ y < x + \eta\},
	\\
	\label{eq:R2}
	&R_2 = \{(x,y)\in\mb R^2, \  \ y > x + \eta\},
\end{align}
where $f_1$ and $f_2$ apply, respectively. 

We now perform a bifurcation analysis of the piecewise-smooth AMOC model~\eqref{model:PWS}. To this end, we use tools from the bifurcation theory for this class of non-smooth systems from the relevant literature \cite{bernardo2008piecewise, filippov2013differential, guardia2011generic}, which we largely follow also in terms of notation and where more details can be found. More specifically, we determine analytic expressions for all (non-smooth) bifurcations, which is possible because of the simple expression for the switching manifold, and the fact that $f_1$ and $f_2$ are linear. We present these results in the form of propositions, whose proofs can be found in Appendix A. The associated curves of codimension-one bifurcations divide the $(\mu,\eta)$-plane into eight open regions, denoted $\mathrm{I-VIII}$. We also present the corresponding phase portraits in the $(\mu,\eta)$-plane, as well as those at the different types of bifurcations. The vertical mixing coefficients are fixed here to $\kappa_1=0.1$ and $\kappa_2=1.0$. This choice is suitable for our purposes and in the realistic range \cite{PCessi1996}, yet slightly different from the values found in the literature \cite{leifeld2016nonsmooth, Walender1982}. Moreover, as can be seen from the expressions in Section~\ref{section:codim1-2}, the bifurcation diagram in the $(\mu,\eta)$-plane is qualitatively the same for any $0<\kappa_1<\kappa_2$.

\subsection{Sliding properties and pseudo-equilibria}
\label{section:fsystem}

We start by introducing some relevant notions from the theory of PWS systems. An equilibrium $p_i$ of the vector field $f_i$ that lies in region $R_i$ is an equilibrium of the overall system and called \emph{admissible}. An important part of the bifurcation theory of planar Filippov systems is the interaction of equilibria and other invariant objects of $f_1$ and $f_2$ with the switching manifold $\Sigma$ \cite{bernardo2008piecewise}. First of all, orbits may cross the switching manifold at the \emph{crossing segment} $\Sigma_c \subset \Sigma$, along which the vector fields $f_1$ and $f_2$ are both transverse and have the same sign. The set of points where $f_1$ and $f_2$ are transversal but have opposite signs is the \emph{sliding segment} $\Sigma_s \subset \Sigma$, which we also refer to as $\Sigma_s^a$ when it is attracting and as $\Sigma_s^r$ when it is repelling. These different segments of the switching manifold are bounded by \emph{tangency points} $F_1$ and $F_2$, where either $f_1$ or $f_2$ is tangent to $\Sigma$, respectively. Generically, such a tangency of $f_i$ is quadratic and isolated, and it is called visible if nearby parabolic orbits lie in $R_i$, and invisible otherwise. For system~\eqref{model:PWS} we have the following.

\begin{prop}[Tangency points and sliding segments]  System~\eqref{model:PWS} has a single sliding segment $\Sigma_s$ that is delimited by two tangency points $F_1$ and $F_2$ at
	\label{prop:escaping_sliding} 
	\begin{align}
		F_i = \binom{1 - \mu + \eta\kappa_i}{1 - \mu + (1 + \kappa_i)\eta} \in \Sigma, 
		\label{pws:tangent_points}
	\end{align}
	which are quadratic when 
	\begin{align}
		\label{cond:Fi_quad}
		\mu + (\mu-\eta-1) \kappa_i - \eta\kappa_i^2 \neq 0.
	\end{align}
	The (quadratic) tangency point $F_1$ is visible for
	\begin{align}
		\label{cond:F1_visible}
		 \mu + (\mu-\eta-1) \kappa_i - \eta\kappa_i^2 < 0,
	\end{align}
	and the (quadratic) tangency point $F_2$ is visible for
	\begin{align}
		\label{cond:F2_visible}
		 \mu + (\mu-\eta-1) \kappa_i - \eta\kappa_i^2 > 0.
	\end{align}
	Otherwise, the (quadratic) tangency at $F_i$ is invisible. For $\eta \neq 0$, system~\eqref{model:PWS} has a sliding segment $\Sigma_s$. When $\eta>0$ the sliding segment is attracting, denoted $\Sigma_s^a$ and given by 
	\begin{align}
		\label{eq:attacking_sliding}
		\Sigma_s^a = \{s \in \Sigma, F_1 < s < F_2\}, 
	\end{align}
	and when $\eta<0$ it is repelling, denoted $\Sigma_s^r$ and given by 
	\begin{align}
		\label{eq:repelling_sliding}
		\Sigma_s^r = \{s \in \Sigma, F_2 < s <F_1\}.
	\end{align}
Here, in a slight abuse of notation, we mean the ordering on the line $\Sigma_s$, as given by the $x$-component.
\end{prop}

A crucial ingredient of the theory is the extension of the flow to the sliding segment $\Sigma_s$ by defining the sliding vector field $f_s$. This is achieved with Filippov's convex method by forming a weighted sum of the adjoining vector fields $f_1$ and $f_2$ such that $f_s$ is in the direction of (the tangent to) $\Sigma_s$  \cite{bernardo2008piecewise, filippov2013differential, guardia2011generic}. With this definition, a PWS orbit is the union of orbit segments induced by the vector fields $f_1$ on $R_1$, $f_2$ on $R_2$, and $f_s$ on $\Sigma_s$. Moreover, every point of the phase plane lies on a unique PWS orbit of the planar Filippov system; see \cite{guardia2011generic} for details. Orbits that remain in $R_1 \cup R_2 \cup \Sigma_c$ are called \emph{regular}, and orbits with segments on $\Sigma_s$ are called \emph{sliding orbits}. Here we use a common convention that sliding orbits continue into $R_1$ or $R_2$ when the end of the sliding segment $\Sigma_s$ is reached (in forward or backward time, respectively, by following the trajectory from the respective tangency point) \cite{guardia2011generic, kuznetsov2003one}. However, the end of $\Sigma_s$ may not be reached because the sliding vector field may have equilibria, called \emph{pseudo-equilibria}, which are referred to as \emph{admissible} when they lie $\Sigma_s$. 

The properties of all equilibria of system~\eqref{model:PWS} can be stated as follows.

\begin{prop}[Equilibria, sliding vector field and pseudo-equilibria]
	\label{prop:equilibria_stabilities}
	System~\eqref{model:PWS} has the following equilibria and pseudo-equilibria for $0<\kappa_1<\kappa_2$. 
	\begin{enumerate}
				\item The vector field $f_i$ has the stable nodal equilibrium
		\begin{align}
			\label{eq:filippov_eq}
			p_i = \left(\frac{1}{1+\kappa_i}, \frac{\mu}{\kappa_i}\right).
		\end{align}
		The equilibrium $p_1$ is admissible when 
		\begin{align}
		\label{cond:p1_admissible}
			\eta>\frac{\mu}{\kappa_1}-\frac{1}{1+\kappa_1},
		\end{align}
		and $p_2$ is admissible when 
		\begin{align}
		\label{cond:p2_admissible}
			\eta<\frac{\mu}{\kappa_2}-\frac{1}{1+\kappa_2}.
		\end{align}
		The admissible equilibrium $p_i$ has a strong stable manifold $W^{ss}(p_i)$ defined by the piecewise-smooth orbit along the linear strong stable direction $W^{ss}_{loc}=\textnormal{span}\binom{1}{0}$. 
		\item The \emph{sliding vector field}  defined on the sliding segment $\Sigma_s$ is given by
\begin{align}
	\label{eq:extended_vector}
        f_s(x)
        &= \frac{1}{\eta} \left(\mu + (\mu - \kappa_2\eta - 1) x + x^2 \right) \binom{1}{1},
\end{align}
where $\Sigma$ is parametrised by $x$.
		\item There are two pseudo-equilibria, that is, equilibria of $f_s$, given by 
		\begin{align}
			\label{eq:pseudo_eq}
			&q^{\pm} = \frac{1}{2}\left(1-\mu-\eta 
\pm\sqrt{(\eta +\mu+1)^2 - 4 \mu} \right) \binom{1}{1} + \binom{0}{\eta}.
		\end{align}
		When $\eta >0$, the pseudo-equilibrium $q^{-}$ is asymptotically unstable and $q^{+}$ is asymptotically stable on $\Sigma_s$. On the other hand, when $\eta<0$, the pseudo-equilibrium $q^{+}$ is asymptotically unstable and $q^{-}$ is asymptotically stable on $\Sigma_s$. The admissibility of these pseudo-equilibria is presented and described in Section~\ref{section:generic_phaseportraits}.
	\end{enumerate}
\end{prop}

Global invariant manifolds of admisible equilibria are defined in complete analogy to those of smooth systems, but with regard to the piecewise-smooth flow $\varphi^t$ constructed in \cite{guardia2011generic}.  Each admissible equilibrium $p_i \in R_i$ of system~\eqref{model:PWS} is attracting with real eigenvalues and, hence, has a strong stable manifold $W^{ss}(p_i)$ consisting of the two orbits that approach $p_i$ tangent to the strong eigenspace. Since system~\eqref{model:PWS} is piecewise linear, $W^{ss}(p_i)$ is actually a straight line locally near $p_i$; however, this is not the case globally since the strong stable manifold typically crosses the switching manifold $\Sigma$. 

We also consider here global invariant manifolds of admissible pseudo-equilibria, which we define as follows. If $q \in \Sigma_s$ is a saddle pseudo-equilibrium then its stable manifold $W^{s}(p)$ or unstable manifold $W^{u}(p)$ is the union of the two \emph{arriving or departing orbits} in $R_1$ and $R_2$, consisting of points that reach $q$ under the piecewise-smooth flow $\varphi^t$ in finite forward or backward time, respectively. The saddle pseudo-equilibrium $q$ then also has associated \emph{generalised (un)stable manifolds} $W^u_{g}(q)$ or $W^s_{g}(q)$. These generalised manifolds consist of segments on $\Sigma_s$ of points that converge to $q$ under the sliding flow (in backward and forward time, respectively), together with their globalisation under $\varphi^t$, which generally consists of departing and arriving orbits to tangency points that bound $\Sigma_s$. When an admissible pseudo-equilibrium $q$ is a nodal attractor, its arriving orbits form the strong stable manifold $W^{ss}(p)$; similarly, a nodal repellor $q \in \Sigma_s$ has the strong unstable manifold $W^{uu}(p)$ consisting of its pair of departing orbits.

\subsection{Bifurcation diagram and structurally stable phase portraits}
\label{section:generic_phaseportraits}

The bifurcation diagram of system~\eqref{model:PWS} consists of curves of (piecewise-smooth) bifurcations that divide the $(\mu, \eta)$-plane into eight open regions $\mathrm{I}$ to $\mathrm{VIII}$, which are equivalence classes of topological equivalence where the phase portraits are structurally stable. This classification is based on the following common notion \cite{guardia2011generic, kuznetsov2003one}: two planar Filippov systems $f$ and $\tilde f$ with switching manifolds $\Sigma$ and $\widetilde \Sigma$, respectively, are \emph{topologically equivalent} if there exists an orientation preserving homeomorphism $h: \mb R^2 \arrow \mb R^2$ that maps $\Sigma$ to $\widetilde \Sigma$ and orbits of $f$ to orbits of $\tilde f$. Note that this definition is a direct and natural extension of that for smooth dynamical systems. In particular, a bifurcation of a planar Filippov system concerns a topological change, and its codimension is given (colloquially speaking) by the number of parameters one needs to find it generically at an isolated point. More information and formal definitions can be found as part of the broad classification in \cite{kuznetsov2003one} of discontinuity-induced bifurcations in planar Filippov systems. 

\begin{figure}[t!]
\centering
	\includegraphics{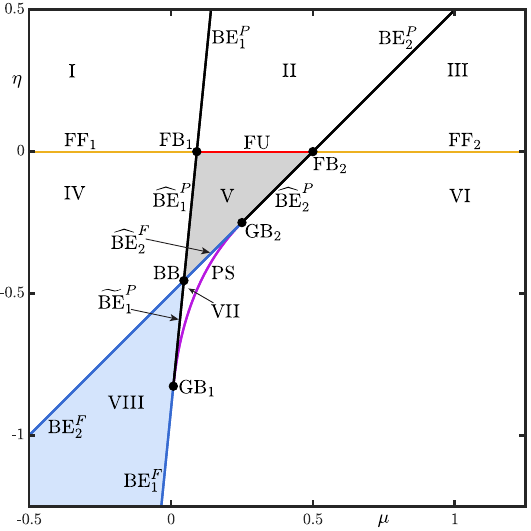}
	\caption{Two-parameter bifurcation diagram of system~\eqref{model:PWS} in the $(\mu, \eta)$-plane with $\kappa_1=0.1$ and $\kappa_2=1.0$. The curves of boundary equilibrium bifurcation $\mathrm{BE}_1$ and $\mathrm{BE}_2$, fold-fold bifurcation $\mathrm{FF}$, and pseudo-saddle-node bifurcation $\mathrm{PS}$ from Proposition~\ref{prop:codim_1} bound regions $\mathrm{I}$ to $\mathrm{VIII}$ with structurally stable phase portraits shown in Figures~\ref{fig:pws_pp1}--\ref{fig:pws_pp3}. Grey shading indicates the existence of a (crossing) periodic orbit, and blue shading bistability between equilibria. These curves intersect at the codimension-two points $\mathrm{FB_1}$, $\mathrm{FB_1}$, $\mathrm{BB}$, $\mathrm{GB_1}$, and $\mathrm{GB_2}$ from Proposition~\ref{prop:pws_codim2}, which generates segments of different bifurcation types shown in Figures~\ref{fig:fold_fold}--\ref{fig:non_smooth_fold}. Specifically, the curves $\mathrm{BE}_1$ and $\mathrm{BE}_2$ consist of segments $\mathrm{BE}_1^P$, $\mathrm{BE}_2^P$, $\widehat{\mathrm{BE}}_1^P$, $\widehat{\mathrm{BE}}_2^P$ and $\widetilde{\mathrm{BE}}_1^P$ of persistence boundary equilibrium bifurcation, and $\mathrm{BE}_1^F$, $\mathrm{BE}_2^F$ and $\widehat{\mathrm{BE}}_2^F$ of non-smooth fold boundary equilibrium bifurcation. The curve $\mathrm{FF}$ consists of segments $\mathrm{FF_1}$ and $\mathrm{FF_2}$ of fold-fold bifurcation and $\mathrm{FU}$ of fused-focus bifurcation.}
	\label{fig:BifurcationDiag}
\end{figure}

Figure~\ref{fig:BifurcationDiag} shows the bifurcation diagram of system~\eqref{model:PWS} in the $(\mu, \eta)$-plane, for the fixed values $\kappa_1=0.1$ and $\kappa_2=1.0$ of the vertical mixing rates, with the 
regions $\mathrm{I}$ to $\mathrm{VIII}$. Their boundaries are formed by bifurcation curves $\mathrm{BE}_1$ and $\mathrm{BE_2}$ of boundary equilibrium bifurcation, $\mathrm{FF}$ of fold-fold bifurcation and $\mathrm{PS}$ of pseudo-saddle-node bifurcation that are formally presented and determined in Proposition~\ref{prop:codim_1}. More precisely, these curves cross or meet at codimension-two bifurcation points $\mathrm{FB_1}$, $\mathrm{FB_1}$, $\mathrm{BB}$, $\mathrm{GB_1}$, and $\mathrm{GB_2}$. As is spelled out in Proposition~\ref{prop:pws_codim2}, these points divide the curves of codimension-one bifurcations into the segments of different bifurcation types that are shown and labeled in Figure~\ref{fig:BifurcationDiag}. 

\begin{figure}[ht!]
\centering
\includegraphics{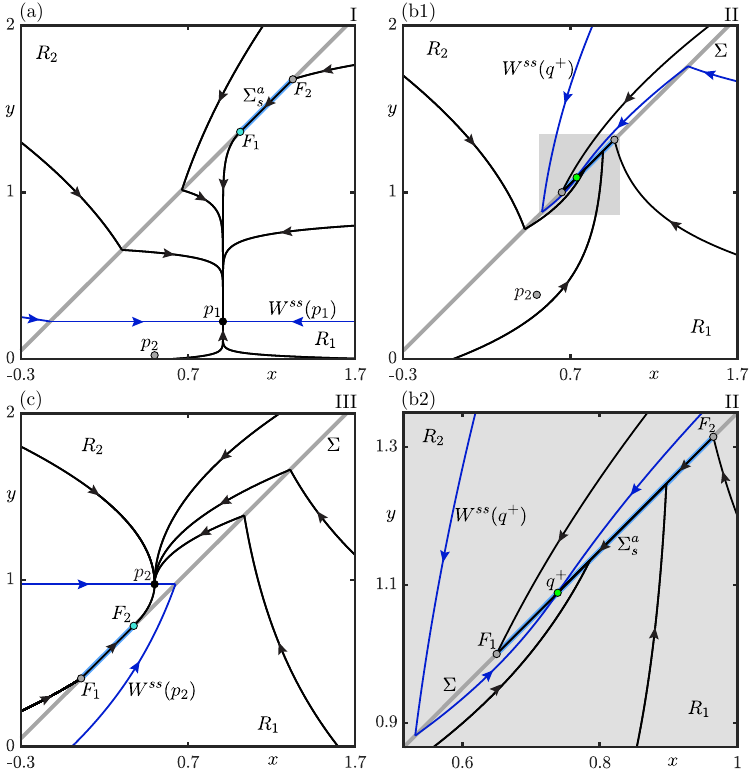}
\caption{Representative phase portraits in regions $\mathrm{I}$ to $\mathrm{III}$ along the horizontal slice $\eta = 0.35$ of the $(\mu, \eta)$-plane, each with an attracting sliding segment $\Sigma_s^a$ bounded by quadratic tangency points $F_1$ and $F_2$. Panel~(a) for $\mu=0.0225$ shows the admissible equilibrium $p_1$ with its strong stable manifold $W^{ss}(p_1)$, as well as the non-admissible equilibrium $p_2 \in R_1$. Panel~(b1) for $\mu = 0.385$ and magnification (b2) near the sliding segment show $p_2 \in R_1$ and the attracting pseudo-node $q^{+}$ with strong stable manifold $W^{ss}(q^{+})$. Panel~(c) for $\mu = 0.975$ shows the admissible equilibrium $p_2 \in R_2$ with $W^{ss}(p_2)$.} 
\label{fig:pws_pp1}	
\end{figure}

\begin{figure}[t!]
	\centering
	\includegraphics{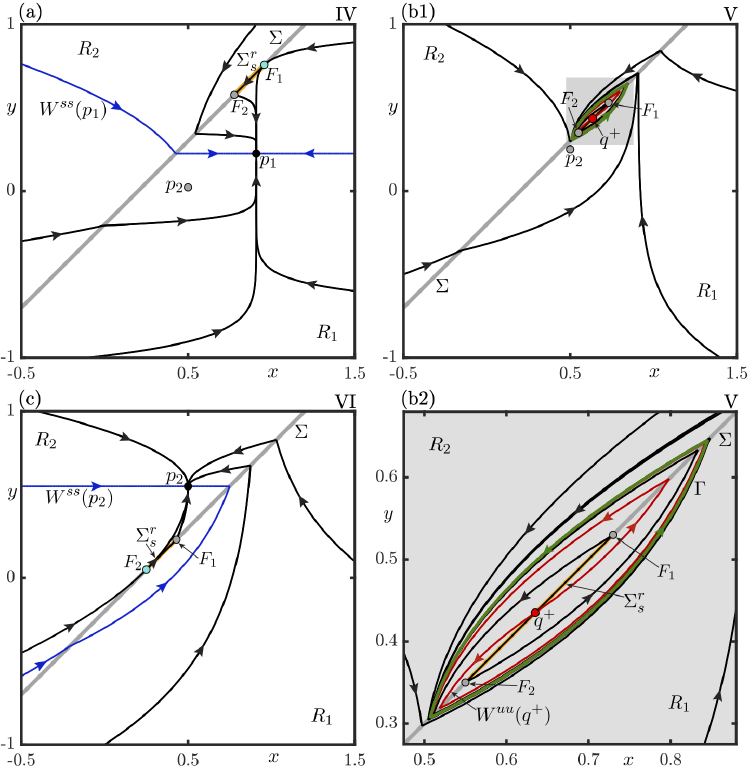}
	\caption{Representative phase portraits in regions $\mathrm{IV}$ to $\mathrm{VI}$ along a horizontal slice $\eta = -0.2$ of the $(\mu, \eta)$-plane, each with a repelling sliding segment $\Sigma_s^r$ bounded by quadratic tangency points $F_1$ and $F_2$. Panel~(a) for $\mu = 0.0225$ is similar to Figure~\ref{fig:pws_pp1}(a), but now the sliding segment is repelling. Panel~(b1) for $\mu=0.25$ and the magnification near the sliding segment (b2) features a repelling pseudo-node $q^{+}$ with a strong unstable manifold $W^{uu}(q^{+})$ and a (crossing) periodic orbit $\Gamma$ that encircles $\Sigma_s^r$. Panel~(c) for $\mu=0.525$ is similar to Figure~\ref{fig:pws_pp1}(c), but now the sliding segment is repelling.}
	\label{fig:pws_pp2}
\end{figure}

\begin{figure}[ht!]
\centering
\includegraphics{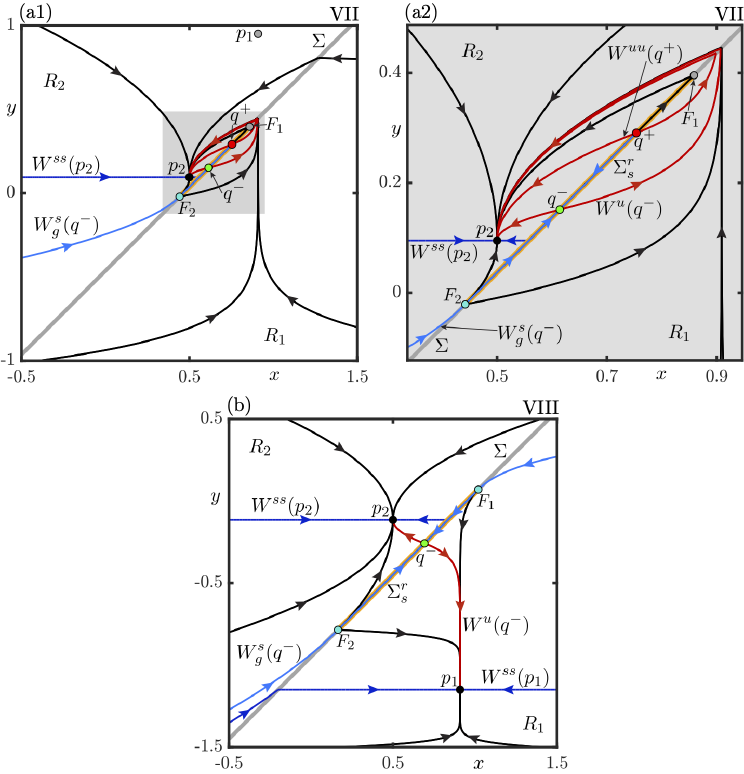}
\caption{Representative phase portraits in regions $\mathrm{VII}$ and $\mathrm{VIII}$, both with a repelling sliding segment $\Sigma_s^r$ bounded by quadratic tangency points $F_1$ and $F_2$. Panel~(a1) for $(\mu,\eta)=(0.0987, -0.463)$ and the magnification near the sliding segment (a2) show $p_2$ with $W^{ss}(p_2)$, pseudo-saddle-equilibrium $q^{-}$ with a generalised stable manifold  $W^{s}_g(q^{-})$ and unstable manifold $W^{u}(q^{-})$, and the repelling pseudo-node $q^{+}$ with strong unstable manifold $W^{uu}(q^{+})$. Panel~(b) for $(\mu,\eta)=(-0.115, -0.95)$ shows the simultaneously admissible equilibria $p_1$ with $W^{ss}(p_1)$ and $p_2$ with $W^{ss}(p_2)$, and the pseudo-saddle-equilibrium $q^{-}$ with $W^{s}_g(q^{-})$ and $W^{u}(q^{-})$.}  
\label{fig:pws_pp3}
\end{figure}

We first present and discuss the structurally stable phase portraits of system~\eqref{model:PWS} in regions $\mathrm{I}$ to $\mathrm{VIII}$; the different bifurcations between them are analysed and illustrated in subsequent sections. The eight cases of phase portraits are shown in Figures~\ref{fig:pws_pp1}--\ref{fig:pws_pp3} in a suitable part of the $(x,y)$-plane. In every phase portrait, the switching manifold $\Sigma$ appears as a straight grey line that partitions phase space into the open regions ${R_1}$ and ${R_2}$. Admissible equilibria of system~\eqref{model:PWS} are shown in black, and non-admissible equilibria in grey. Non-admissible equilibria outside the frame of interest (far away from the switching manifold) are not shown. There exists a sliding segment in each region: attracting sliding segments $\Sigma_s^a$ are coloured blue, and repelling sliding segments $\Sigma_s^r$ are coloured orange. In either case, the sliding segment is bounded by the quadratic tangency points $F_1$ and $F_2$, which are coloured cyan when visible and grey when invisible. Admissible pseudo-equilibria $q^{-}$ and $q^{+}$ are coloured by their stability: stable pseudo-equilibria are green, and unstable pseudo-equilibria are red. Admissible equilibria and pseudo-equilibria may have (strong) invariant manifolds that are coloured blue when stable and red when unstable. Some representative trajectories are shown in black, and they were obtained numerically with an integrator based on event-detection, as described in \cite{piiroinen2008event}. 

Figure~\ref{fig:pws_pp1} presents phase portraits of system~\eqref{model:PWS} in regions $\mathrm{I}$ to $\mathrm{III}$, which all feature and attracting sliding segment $\Sigma_s^a$. In the phase portrait in region $\mathrm{I}$, shown in panel~(a), $\Sigma_s^a$ is bounded by a visible quadratic tangency point $F_1$ on the left and an invisible quadratic tangency point $F_2$ on the right. Neither of the pseudo-equilibria $q^{-}$ and $q^{+}$ are on $\Sigma_s^a$ and, hence, they are non-admissible (and not shown). The equilibrium $p_2$ lies in region $R_1$ and is non-admissible, while $p_1 \in R_1$ is admissible and a global attractor. Orbits in $R_1$ and $R_2$ are either regular and converge to $p_1$ or hit the attracting sliding segment $\Sigma_s^a$, along which sliding orbits approach $F_1$ and then depart into $R_1$ to converge to $p_1$. Note, that the strong stable manifold $W^{ss}(p_1)$ of $p_1$ is composed of a horizontal component in ${R_1}$ and the corresponding arriving orbit in $R_2$. Crossing the segment $\mathrm{BE}_1^P$ of boundary equilibrium bifurcation results in $p_1$ becoming non-admissible by moving into $R_2$ through $F_1$; at the same time, a pseudo-equilibrium $q^{+} \in \Sigma_s^a$ emerges from $F_1$, where the tangency is now invisible. The resulting phase portrait in region $\mathrm{II}$ is shown in Figure~\ref{fig:pws_pp1}(b1) with a magnification near the sliding segment in panel~(b2). The pseudo-equilibrium $q^{+}$ is a global attractor: all orbits in $R_1$ and $R_2$ hit the sliding segment $\Sigma_s^a$, along which the sliding orbits converge to $q^{+}$. Moreover, $q^{+}$ has the strong stable manifold $W^{ss}(q^{+})$, consisting of the two arriving orbits to $q^{+}$ from within $R_1$ and $R_2$, respectively; see panel~(b2). When the segment $\mathrm{BE}_2^P$ is crossed there is again a boundary equilibrium bifurcation, but now of $p_2$ at $F_2$: as Figure~\ref{fig:pws_pp1}(c) shows, in region $\mathrm{III}$ the pseudo-equilibrium $q^{+}$ moved off $\Sigma_s^a$ through $F_2$, and $p_2$ with strong stable manifold $W^{ss}(p_2)$ is now admissible and the global attractor. 

Phase portraits in regions $\mathrm{IV}$ to $\mathrm{VI}$ are presented Figure~\ref{fig:pws_pp2}; as was the case for regions $\mathrm{I}$ to $\mathrm{III}$, this also concerns the the transition from $p_1$ to $p_2$ being the global attractor, with the difference that there is now a repelling sliding segment $\Sigma^r_s$. In the phase portrait in region $\mathrm{IV}$, shown in panel~(a), $\Sigma^r_s$ is bounded by an invisible quadratic tangency point $F_2$ on the left and by a visible quadratic tangency point $F_1$ on the right; the pseudo-equilibria $q^{-}$ and $q^{+}$ are on $\Sigma_c$ and non-admissible (and not shown). The only admissible equilibrium is $p_1 \in R_1$, and it is a global attractor. Sliding orbits on $\Sigma_s^r$ approach $F_2$, where they depart into $R_1$ and converge to $p_1$; however, in contrast to region $\mathrm{III}$, no forward orbits hit $\Sigma_s^r$ as this sliding segment is repelling. When crossing segment $\widehat{\mathrm{BE}}_1^P$, we find again a (persistence) boundary equilibrium bifurcation where $p_1$ moves through $F_1$ and becomes non-admissible. As Figure~\ref{fig:pws_pp2}(b1) and the magnification in panel~(b2) show, in region $\mathrm{V}$ this results again in the pseudo-equilibrium $q^{+}$ being admissible. However, $q^{+} \in \Sigma_s^r$ is now a repelling node with strong unstable manifold $W^{uu}(q^{+})$, consisting of the departing orbits from $q^{+}$ in $R_1$ and $R_2$, respectively. Importantly, in region $\mathrm{V}$ there is a stable (crossing) periodic orbit $\Gamma$, which is composed of orbit segments of $f_1$ in $R_1$ and $f_2$ in $R_2$ that join on the crossing segment $\Sigma_c$. All points except $q^{+}$ converge to this periodic orbit; in particular, $W^{uu}(q^{+})$ accumulates on $\Gamma$, while initial conditions on $\Sigma_s^r \setminus \{q^{+}\}$ move to an end point $F_2$ or $F_1$ of $\Sigma_s^r$, where they depart into $R_1$ or $R_2$, respectively, to converge to $\Gamma$; see panel~(b2). Crossing segment $\widehat{\mathrm{BE}}_2^P$ concerns a second (persistence) boundary equilibrium bifurcation, but now of $p_2$ at the tangent point $F_2$. As a result, the now admissible equilibrium $p_2$ is indeed the global attractor in region $\mathrm{VI}$, as is shown in Figure~\ref{fig:pws_pp2}(c). 

Phase portraits for regions $\mathrm{VII}$ and $\mathrm{VIII}$ are presented in Figure~\ref{fig:pws_pp3}; they both still feature the repelling sliding segment $\Sigma^r_s$, bounded by the quadratic tangency points $F_2$ on the left and $F_1$ on the right. In region $\mathrm{VII}$, as in panel~(a1) with a magnification in panel~(a2), we find the repelling pseudo-equilibrium $q^{+} \in \Sigma_s^r$ as in region $\mathrm{V}$. However, due to the transition through the bounding segment $\widehat{\mathrm{BE}}_2^F$ of (non-smooth fold) boundary equilibrium bifurcation, the equilibrium $p_2$ is now in open region $R_2$ and admissible, and the second pseudo-equilibrium $q^{-}$ now also lies on $\Sigma_s^r$. The point $p_2$ attracts all points, apart from those on the generalised stable manifold $W^{s}_g(q^{-})$ of $q^{-}$, which is composed of sliding orbits on $\Sigma^r_s$ approaching $q^{-}$ and the arriving orbit to $F_2$. The unstable manifold $W^u(q^{-})$ of $q^{-}$ and strong unstable manifold $W^{uu}(q^{+})$ both converge to the attractor $p_2$. Note that $W^{ss}(p_2)$ is composed of a horizontal component in $R_1$, the corresponding sliding orbit in $\Sigma^r_s$ and the arriving orbit to $F_1$. When crossing segment $\widetilde{\mathrm{BE}}_1^P$ into region $\mathrm{VIII}$, there is a (persistence) boundary equilibrium bifurcation, at which $p_1$ becomes admissible and the pseudo-equilibrium $q^{+}$ becomes non-admissible by moving through $F_1$ onto the crossing segment $\Sigma_c$. As the phase portrait in Figure~\ref{fig:pws_pp3}(b) shows, both $p_1 \in R_1$ and equilibrium $p_2 \in R_2$ are now attractors in region $\mathrm{VIII}$; hence, this is the region of bistability. The generalised stable manifold $W_g^{s}(q^{-})$ of the saddle pseudo-equilibrium $q^{-} \in \Sigma^r_s$ is now composed of $\Sigma^r_s$ and the arriving orbits to both $F_1$ and $F_2$, and it forms the boundary between the basins of attraction of the attractors $p_1$ and $p_2$.  Indeed, the lower branch of $W^{u}(q^{-})$ converges to $p_1$, and its upper branch to $p_2$.

\subsection{Codimension-one and codimension-two bifurcations}
\label{section:codim1-2}

We now present analytical expressions for all (non-smooth) bifurcations of system~\eqref{model:PWS} of codimension one and two in Propositions~\ref{prop:codim_1} and~\ref{prop:pws_codim2}, respectively. The respective proofs can be found in Appendix A.

\begin{prop}[Codimension-one bifurcations]
	\label{prop:codim_1}
	System~\eqref{model:PWS} has the following codimension-one bifurcations, further information on which can be found in \cite{kuznetsov2003one}. 
		\begin{enumerate}
		\item \emph{Boundary equilibrium bifurcations} of $p_1$ and $p_2$ occur, respectively, along the straight lines 
		\begin{align*}
			&\mathrm{BE}_1: \ \ (\mu,\eta) = \left(\mu, \ \frac{\mu}{\kappa_1} - \frac{1}{1+\kappa_1}\right),
			\\
			&\mathrm{BE}_2: \ \ (\mu,\eta) = \left(\mu, \ \frac{\mu}{\kappa_2} - \frac{1}{1+\kappa_2}\right).
		\end{align*}
		At $\mathrm{BE}_i$ the equilibrium $p_i$ collides with the tangency point $F_i$, changing its visibility. The tangency point $F_1$ is visible in regions $\mathrm{I, IV}$ and $\mathrm{VIII}$. Similarly, the tangency point $F_2$ is visible in regions $\mathrm{III, VI, VII}$ and $\mathrm{VIII}$. The pseudo-equilibrium $q^{+}$ is admissible in regions $\mathrm{II}$, $\mathrm{V}$ and $\mathrm{VII}$. Similarly, the pseudo-equilibrium $q^{-}$ is admissible in regions $\mathrm{VII}$ and $\mathrm{VIII}$.
		\item \emph{Fold-fold bifurcations} occur along the horizontal line 
		\begin{align*}
			\mathrm{FF}: \ \ (\mu,\eta) = (\mu, \ 0).
		\end{align*}
		At $\mathrm{FF}$ the tangency points $F_1$ and $F_2$ coincide at a singular tangency point $F^*$ and switch places on the sliding segment boundary, resulting in the sliding segments changing between being attracting and repelling \cite{kuznetsov2003one}; see also Proposition~\ref{prop:escaping_sliding} for a description of the tangency points.
		\item  A \emph{pseudo-saddle-node bifurcation} occurs along the curve segment
		\begin{align*}
			&\mathrm{PS}: \ \ (\mu,\eta)=(\mu, \ -(\mu+1) + 2\sqrt{\mu}), \ \ \ \frac{\kappa_1^2}{(\kappa_1+1)^2} < \mu < \frac{\kappa_2^2}{(\kappa_2^2+1)^2}.
		\end{align*}
		Along $\mathrm{PS}$ the pseudo-equilibria $q^{-}$ and $q^{+}$ form a saddle-node at 
		\begin{align*}
\label{eq:qstar}
			& q^*=(1-\sqrt{\mu}) \binom{1}{1} + \binom{0}{\eta}
		\end{align*}
 on the repelling sliding segment $\Sigma_s^r$. 
	\end{enumerate}
\end{prop}

The curves $\mathrm{BE}_1$, $\mathrm{BE}_2$, $\mathrm{FF}$ and $\mathrm{PS}$ from Proposition~\ref{prop:codim_1} intersect or meet at codimension-two bifurcation points. These points divide $\mathrm{BE}_1$, $\mathrm{BE}_2$, $\mathrm{FF}$ into the segments shown in Figure~\ref{fig:BifurcationDiag}, along which the respective codimension-one bifurcation manifests itself in a topologically different way, as follows.

\begin{prop}[Codimension-two bifurcations] 
	\label{prop:pws_codim2}
	System~\eqref{model:PWS} has the following codimension-two bifurcations for $0<\kappa_1<\kappa_2$. 
	\begin{enumerate}
		\item \emph{Fold-boundary equilibrium bifurcations}
		\begin{align}
			&\mathrm{FB_1}: \ \ (\mu,\eta) = 
			\left(  
			\frac{\kappa_1}{1 + \kappa_1},\ 0
			\right),
			\\
			&\mathrm{FB_2}: \ \ (\mu,\eta) = 
			\left(  
			\frac{\kappa_2}{1 + \kappa_2},\ 0
			\right),
		\end{align}
		occur at the intersection point of the curve $\mathrm{FF}$ with the curves $\mathrm{BE_1}$ and $\mathrm{BE_2}$, respectively. At the point $\mathrm{FB}_i$, the equilibrium $p_i$ collides with the singular tangency point $F^*$. 
		\\
		\\
		The point $\mathrm{FB}_i$ divides the curve $\mathrm{FF}$ locally into segments $\mathrm{FF}_1$ and $\mathrm{FF_2}$, which is the case of fold-fold bifurcation of type VI$_1$ as presented in \cite{kuznetsov2003one}, and the segment $\mathrm{FU}$ of fused-focus bifurcation \cite{kuznetsov2003one, castillo2017pseudo} along which the (crossing) periodic orbit $\Gamma$ (dis)appears. Both of these fold-fold bifurcations result in the sliding segment changing between being repelling and attracting, and the quadratic tangency points $F_i$ switching places as the sliding segment boundaries. 
		\\
		\\
		The point $\mathrm{FB}_i$ also divides the curve $\mathrm{BE}_i$ locally into segment $\mathrm{BE}^P_i$, where there is a standard persistence boundary equilibrium bifurcation with a nodal equilibrium as presented in \cite{kuznetsov2003one}, and a segment  $\widehat{\mathrm{BE}}_i^P$ along which a stable (crossing) periodic orbit $\Gamma$ (dis)appears in a homoclinic-like persistence boundary equilibrium bifurcation.
		\\
		
		\item A \emph{double-boundary equilibrium bifurcation}
		\begin{align}
			&\mathrm{BB}: \ \ (\mu,\eta) = 
			\left( \frac{\kappa_1\kappa_2}{(\kappa_1 + 1)(\kappa_2 + 1)},\ \ \frac{-1}{(\kappa_1 + 1)(\kappa_2 + 1)}
			\right),
		\end{align}
		occurs at the intersection of the curves $\mathrm{BE_1}$ and $\mathrm{BE_2}$. At point $\mathrm{BB}$ the equilibria $p_1$, and $p_2$ simultaneously collide at the two different quadratic tangency points $F_1$ and $F_2$, respectively. 
		\\
		\\
		The point $\mathrm{BB}$ divides the curve $\mathrm{BE_1}$ locally into segment $\widehat{\mathrm{BE}}_1^P$, along which a stable (crossing) periodic orbit $\Gamma$ (dis)appears, and a segment $\widetilde{\mathrm{BE}}_1^P$, where there is a standard persistence boundary equilibrium bifurcation with a nodal equilibrium as presented in \cite{kuznetsov2003one}. Similarly, the curve $\mathrm{BE_2}$ is divided locally by $\mathrm{BB}$ into the segment $\widehat{\mathrm{BE}}_2^F$, along which the (crossing) periodic orbit $\Gamma$ (dis)appears, and a segment $\mathrm{BE}^F_2$, where there is the standard non-smooth fold boundary equilibrium bifurcation with a nodal equilibrium  \cite{kuznetsov2003one}.
		\\
		
		\item \emph{Generalized boundary equilibrium bifurcations} \cite{guardia2011generic} 
		\begin{align}
			&\mathrm{GB_1}: \ \ (\mu, \eta) = 
			\left(  
			\frac{\kappa_1^2}{(\kappa_1 + 1)^2}, \ \  -\frac{1}{(\kappa_1+1)^2}
			\right),
			\\
			&\mathrm{GB_2}: \ \ (\mu, \eta) = 
			\left(
			\frac{\kappa_2^2}{(\kappa_2 + 1)^2}, \ \ -\frac{1}{(\kappa_2+1)^2}
			\right),
		\end{align}
		occur at end points of the curve $\mathrm{PS}$, respectively, on the curves $\mathrm{BE}_1$ and $\mathrm{BE}_2$. At the point $\mathrm{GB}_i$, equilibrium $p_i$ collides with the quadratic tangency point $F_i$. At the same time, a pseudo-saddle-node bifurcation takes place at $F_i$, resulting in a generalised boundary equilibrium bifurcation with respect to $f_i$. 
		\\ 
		\\
		The point $\mathrm{GB_1}$ separates the curve $\mathrm{BE_1}$ locally into the segment $\widetilde{\mathrm{BE}}^P_1$ and the segment ${\mathrm{BE}}^F_1$ of boundary equilibrium bifurcations. Similarly, the point $\mathrm{GB_2}$ separates the curve $\mathrm{BE_2}$ locally into the segment $\widehat{\mathrm{BE}}^P_2$ and the segment $\widehat{\mathrm{BE}}^F_2$ of boundary equilibrium bifurcations. 
	\end{enumerate}
\end{prop}

\subsection{Phase portraits at codimension-one bifurcations}
\label{section:pp_codim1}

We now present in Figures~\ref{fig:fold_fold}--\ref{fig:non_smooth_fold} phase portraits in the $(x,y)$-plane for each segment of codimension-one bifurcation introduced in Proposition~\ref{prop:pws_codim2} and shown and labeled accordingly in Figure~\ref{fig:BifurcationDiag}. Here, we take a global view of each such transition to allow for comparison with the respective neighbouring structurally stable phase portraits in Figures~\ref{fig:pws_pp1}--\ref{fig:pws_pp3}.

\subsubsection{Fold-fold and pseudo-Hopf bifurcations}

\begin{figure}[t!]
	\centering
	\includegraphics{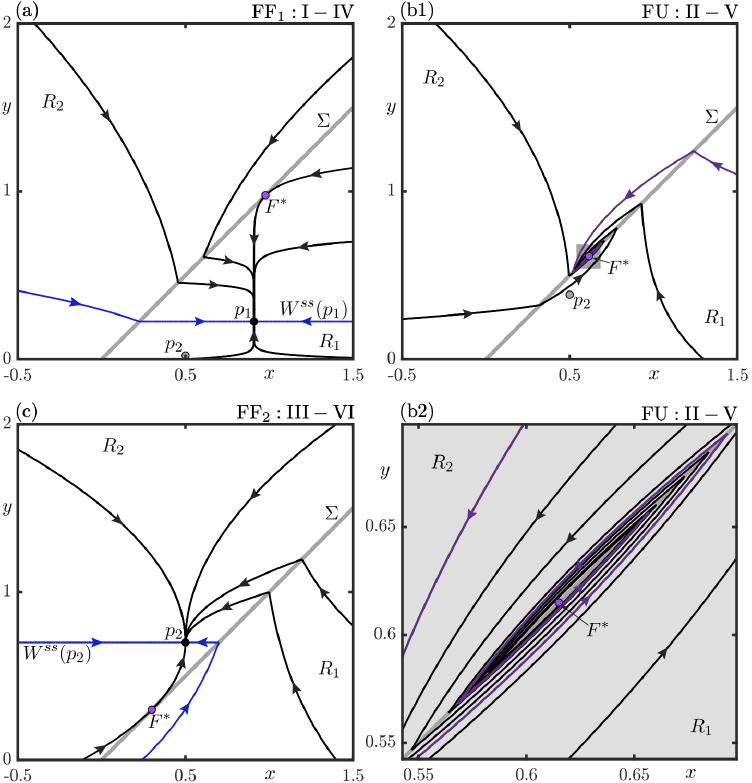}
	\caption{Representative phase portraits of segments $\mathrm{FF_1}, \mathrm{FU}$ and $\mathrm{FF_2}$ along $\eta = 0$, each with a singular tangency point $F^{*}$. Panel~(a) for $\mu = 0.0225$ on $\mathrm{FF}_1$ shows the global attractor $p_1$ with $W^{ss}(p_1)$ and the equilibrium $p_2$. Panel~(b1) for $\mu=0.25$ on $\mathrm{FU}$ and the magnification (b2) shows the weakly attracting fold-fold point $F^*$; a representative orbit is highlighted in purple. Panel~(c) for $\mu=0.7$ on $\mathrm{FF}_2$ shows the global attractor $p_2$ with $W^{ss}(p_2)$.} 
	\label{fig:fold_fold}
\end{figure}

Figure~\ref{fig:fold_fold} shows the phase portraits along segments $\mathrm{FF_1}$, $\mathrm{FU}$, and $\mathrm{FF_2}$, each of which with a singular tangency point $F^*$. The phase portrait along segment $\mathrm{FF_1}$, which separates regions $\mathrm{I}$ and $\mathrm{IV}$, is shown in panel~(a). The admissible equilibrium $p_1 \in R_1$ is a global attractor with strong stable manifold $W^{ss}(p_1)$. The singular tangencyat the point $F^{*}$ is invisible to $f_2$ and visible to $f_1$, and orbits of system~\eqref{model:PWS} are collinear at $F^{*}$. The phase portrait along segment $\mathrm{FU}$, separating regions $\mathrm{II}$ and $\mathrm{V}$, is shown in panel~(b1) with a magnification near $F^{*}$ in panel~(b2). The singular tangency at $F^{*}$ is now invisible to both vector fields $f_1$ and $f_2$, and orbits are anti-collinear at $F^{*}$. Therefore, orbits spiral inward toward $F^{*}$ (at a very slow rate), and this point is a global attractor. This situation is reminiscent of a (supercritical) Hopf bifurcation for smooth dynamical systems, which is why this bifurcation is also known as a pseudo-Hopf bifurcation \cite{castillo2017pseudo}. The phase portrait along segment $\mathrm{FF_2}$, which separates regions $\mathrm{III}$ and $\mathrm{VI}$, is presented in panel~(c). The singular tangency at $F^{*}$ is now visible to $f_2$ and invisible to $f_1$, and $p_2 \in R_2$ is admissible and the global attractor.

\subsubsection{Boundary equilibrium and pseudo-saddle-node bifurcations}

The phase portraits along the segments of the curves $\mathrm{BE}_1$ and $\mathrm{BE}_2$ from Proposition~\ref{prop:codim_1} are characterised by an equilibrium of system~\eqref{model:PWS} being on the switching manifold, but they have different global manifestations.

\begin{figure}[t!]
	\centering
	\includegraphics{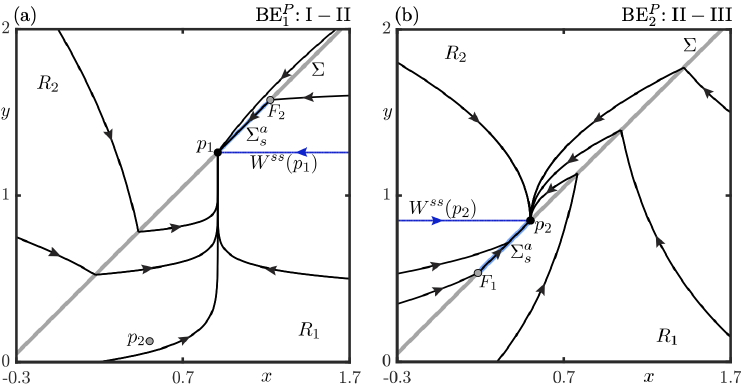}
	\caption{Representative phase portraits along the segments $\mathrm{BE}_1^P$ and $\mathrm{BE}_2^P$ on the horizontal slice $\eta = 0.35$ of the $(\mu,\eta)$-plane, each with an attracting sliding segment $\Sigma_s^a$. Panel~(a) for $\mu = 0.1259$ on $\mathrm{BE}_1^P$ shows boundary equilibrium $p_1$ with $W^{ss}(p_1)$. Panel~(b) for $\mu = -0.25$ on $\mathrm{BE}_2^P$ shows the boundary equilibrium $p_2$ with $W^{ss}(p_2)$.}
	\label{fig:persistance_upper}
\end{figure}

The phase portrait along segment $\mathrm{BE}_1^P$, which separates regions $\mathrm{I}$ via $\mathrm{II}$, is shown in Figure~\ref{fig:persistance_upper}(a). It has the attracting sliding segment $\Sigma_s^a$ bounded by the boundary-node $p_1$ on the left, and by the invisible quadratic tangency point $F_2$ on the right; both pseudo-equilibria $q^{-}$ and $q^{+}$ are non-admissible (and not shown). The equilibrium $p_2$ is not admissible and the boundary-node $p_1$ is a global attractor; note that its strong stable manifold $W^{ss}(p_1)$ consists only of the horizontal arriving orbit to $p_1$ in $R_1$.  The phase portrait in panel~(b) along segment $\mathrm{BE}_2^P$, separating regions $\mathrm{II}$ and $\mathrm{III}$, is the corresponding situation but for the boundary-node $p_2$: this point is now the  global attractor with strong stable manifold $W^{ss}(p_2)$ in $R_2$, and it bounds $\Sigma_s^a$ together with the invisible quadratic tangency point $F_1$. 

\begin{figure}[t!]
\centering
\includegraphics{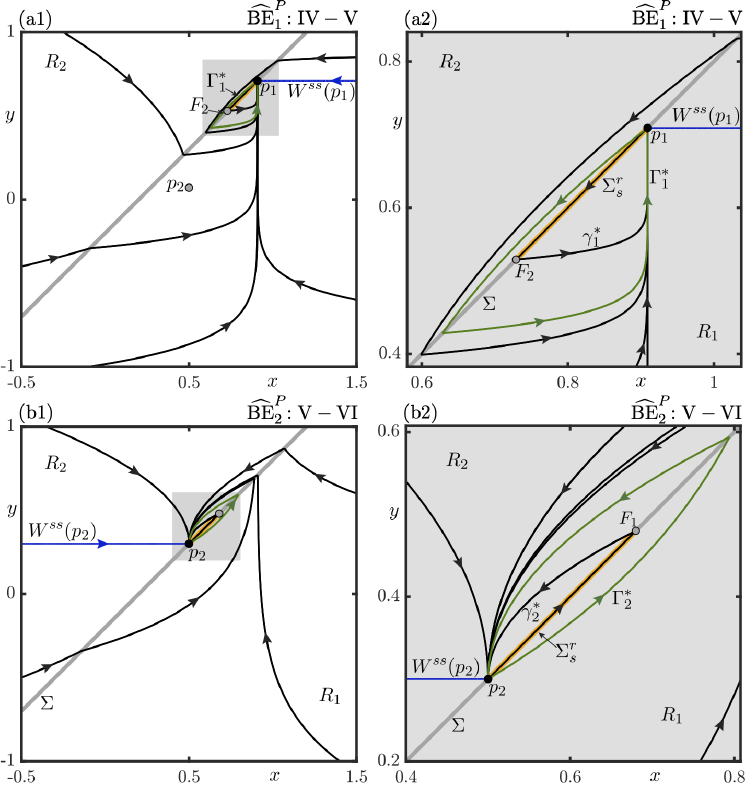}
\caption{Representative phase portraits along segments $\widehat{\mathrm{BE}}_1^P$ and $\widehat{\mathrm{BE}}_2^P$ on the horizontal slice $\eta = -0.2$ in the $(\mu,\eta)$-plane, each with a repelling sliding segment $\Sigma_s^r$. Panel~(a1) for $\mu = 0.0709$ on $\widehat{\mathrm{BE}}_1^P$ and the magnification (a2) show the boundary equilibrium $p_1$ with $W^{ss}(p_1)$, and homoclinic connections $\gamma^*_1$ and $\Gamma^*_1$.  Panel~(b1) for $\mu = 0.3$ on $\widehat{\mathrm{BE}}_2^P$ and the magnification (b2) show boundary equilibrium $p_2$ with $W^{ss}(p_2)$, and homoclinic connections $\gamma_2^*$ and $\Gamma_2^*$.} 
\label{fig:persistance_1}
\end{figure}

Figure~\ref{fig:persistance_1}(a1) shows the phase portrait along segment $\widehat{\mathrm{BE}}_1^P$, separating regions $\mathrm{IV}$ and $\mathrm{V}$, with a magnification in panel~(a2) near the sliding segment, which is now repelling. Here $\Sigma_s^r$ is bounded by the invisible quadratic tangency point $F_2$ on the left and by the boundary-node $p_1$ on the right, with both pseudo-equilibria non-admissible (and not shown). The point $p_1$ is globally attracting with strong stable manifold $W^{ss}(p_1)$ in $R_1$. However, the vector field $f_2$ is transverse to $\Sigma$ at $p_1$, and this departing orbit forms a non-sliding homoclinic connection $\Gamma_1^*$ back to the boundary-node $p_1$. Observe in panel~(a2) that $\Gamma_1^*$ bounds a region of a family of homoclinic orbits that involve sliding (in backward time) along the repelling sliding segment $\Sigma_s^r$. The orbit labeled $\gamma^*_1$, consisting of $\Sigma_s^r$ and the departing orbit from $F_2$, is the maximal sliding homoclinic orbit: it divides this region inside $\Gamma_1^*$ into homoclinic orbit that remain in $R_1$ from those that have segments in both $R_1$ and $R_2$. The phase portrait along segment $\widehat{\mathrm{BE}}_2^P$, which separates regions $\mathrm{V}$ and $\mathrm{VI}$, is shown similarly in Figure~\ref{fig:persistance_1}(b1) and~(b2). The overall picture is effectively that same, but now $p_2$ is the globally attracting boundary equilibrium on $\Sigma_s^r$, with analogous non-sliding and maximal sliding homoclinic orbits $\Gamma_2^*$ and $\gamma^*_2$, respectively. The characterising feature of this bifurcation is the existence of a non-sliding and crossing homoclinic connection $\Gamma_i^*$, from which the stable (crossing) periodic orbit $\Gamma$ in region $\mathrm{V}$ bifurcates; compare with  Figure~\ref{fig:pws_pp2}(b). This type of (persistence) boundary equilibrium bifurcation is hardly discussed in the literature; to our knowledge, it has only been observed in the related Welander's box model in \cite{leifeld2016nonsmooth}, where it is referred to as a homoclinic-like boundary equilibrium bifurcation.  

\begin{figure}[h!]
	\centering
	\includegraphics{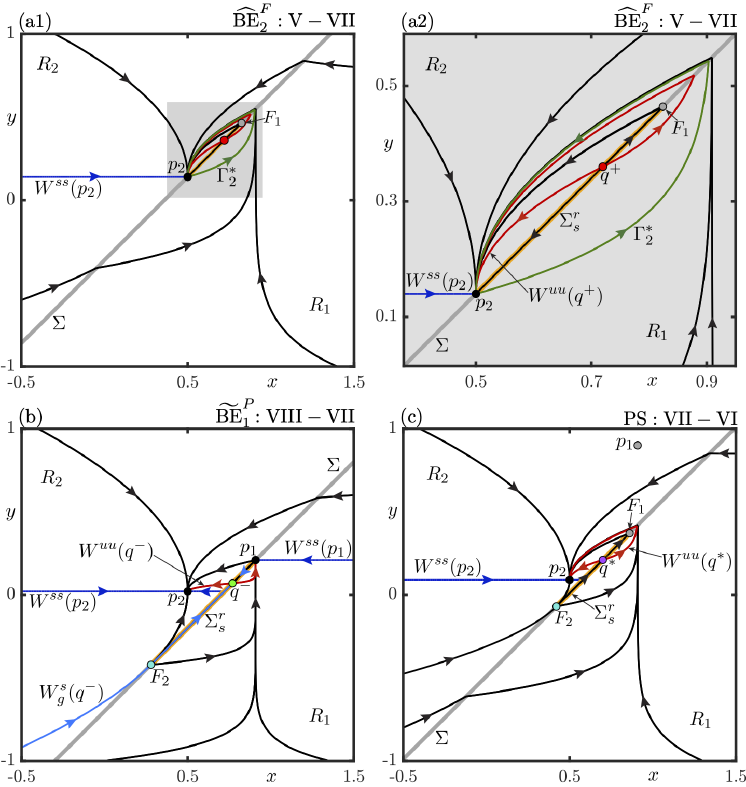}
	\caption{Representative phase portraits along segments $\widehat{\mathrm{BE}}_2^F$, $\widetilde{\mathrm{BE}}_1^P$ and along the curve $\mathrm{PS}$, each with a repelling sliding $\Sigma_s^a$. Panel~(a1) for $(\mu,\eta) = (0.14, -0.36)$ on $\widehat{\mathrm{BE}}_2^F$ and the magnification (a2) show the boundary equilibrium $p_2$ with $W^{ss}(p_2)$ and homoclinic connection $\Gamma_2^*$, and the pseudo-equilibrium $q^{+}$ with $W^{uu}(q^{+})$. Panel~(b) for $(\mu,\eta)=(-0.0141, -1.05)$ on $\widetilde{\mathrm{BE}}_1^P$ shows the admissible equilibrium $p_2 \in R_2$ with $W^{ss}(p_2)$, the boundary equilibrium $p_1$ with $W^{ss}(p_1)$, and the pseudo-equilibrium $q^{-}$ with $W^{u}(q^{-})$ and $W^{s}_g(q^{-})$. Panel~(c) for $(\mu,\eta)=(0.09, -0.7636)$ on $\mathrm{PS}$ shows the equilibrium $p_2$ with $W^{ss}(p_2)$ and the singular pseudo-equilibrium $q^{*}$ with $W^{uu}(q^{*})$.}
	\label{fig:pseudo_saddle_node}
\end{figure}

The phase portrait along segment $\widehat{\mathrm{BE}}^F_2$, which separates regions $\mathrm{V}$ and  $\mathrm{VII}$, is shown in  Figure~\ref{fig:pseudo_saddle_node}(a1) with a magnification near the repelling sliding segment in panel~(a2). Here, $\Sigma_s^r$ is bounded by the attracting boundary-node $p_2$ on the left and by an invisible tangency $F_1$ on the right; moreover, it contains the admissible and repelling pseudo-equilibrium $q^{+} \in \Sigma_s^r$ (while the pseudo-equilibrium $q^{-}$ is non-admissible and not shown). As was the case along segment $\widehat{\mathrm{BE}}_1^F$, the phase portrait in Figure~\ref{fig:pseudo_saddle_node}(a) features a (crossing) homoclinic orbit $\Gamma^*_2$ of $p_2$. However, due to the existence of $q^{+}$ on $\Sigma_s^r$, this special orbit does now not bound a region with further (sliding) homoclinic orbits. Regardless, $\Gamma^*_2$ is still the limit of the stable (crossing) periodic orbit $\Gamma$ in region $\mathrm{V}$. Note that all points inside the region bounded by $\Gamma^*_2$ converge in backward time to the unstable pseudo-equilibrium $q^{+}$, whose strong unstable manifold $W^{uu}(q^{+})$ converges to $p_2$; see panel~(a2). Segment $\widetilde{\mathrm{BE}}_1^P$ separates regions $\mathrm{VII}$ and $\mathrm{VIII}$, and the phase portrait along it is shown in Figure~\ref{fig:pseudo_saddle_node}(b). Here, $p_1$ is the attracting boundary-node, and the quadratic tangency $F_2$ is visible. The equilibrium $p_2 \in R_2$ is admissible and also attracting. Moreover, the pseudo-equilibrium $q^{-}$ lies on the repelling sliding section $\Sigma_s^r$, and it is a saddle. Its generalised stable manifold $W^{s}_g(q^{-})$ consists of $\Sigma_s^r$ and the arriving orbit to $F_2$. The boundary between the basins of attraction of $p_1$ and $p_2$ is formed by the union of $W^{s}_g(q^{-})$ and the strong stable manifold $W^{ss}(p_1)$ in $R_1$. Points below these curves and including $W^{ss}(p_1)$ converge to $p_1$, while points above these curves converge to $p_2$. 

The pseudo-saddle-node bifurcation along the curve $\mathrm{PS}$ separates regions $\mathrm{VI}$ and $\mathrm{VII}$, and its phase portrait is shown in  Figure~\ref{fig:pseudo_saddle_node}(c). As the name suggests, there is a saddle-node  $q^{*}$ of pseudo-equilibria on the repelling sliding segment $\Sigma_s^r$, which is the limiting point where the admissible pseudo-equilibria $q^{-}, q^{+} \in \Sigma_s^r$ in region $\mathrm{VII}$ (dis)appear. Note that $q^{*}$ is semi-stable on $\Sigma_s^r$ and has the strong unstable manifold $W^{uu}(q^*)$. The points on $\Sigma_s^r$ in between the visible quadratic tangency point $F_2$ and $q^{*}$ end up at $q^{*}$ under the sliding flow; all other points in the $(x,y)$-plane converge to the admissible and stable equilibrium $p_2 \in R_2$ with strong stable manifold $W^{ss}(p_2)$.

\begin{figure}[t!]
	\centering
	\includegraphics{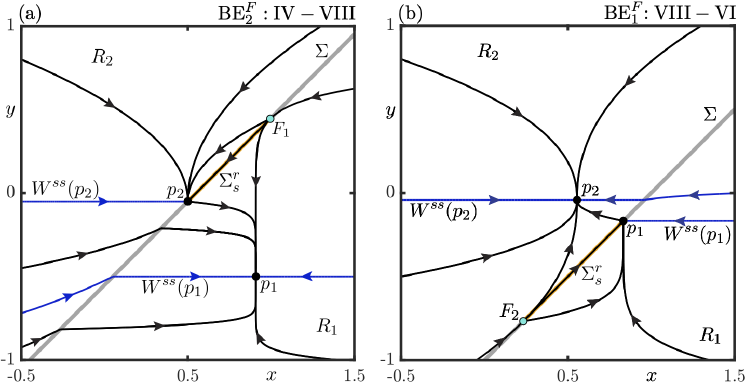}
	\caption{Representative phase portraits along the segments ${\mathrm{BE}}_2^F$ and $\mathrm{BE}_1^F$, each with a repelling sliding segment $\Sigma_s^r$.  Panel~(a) for $(\mu, \eta)=(-0.05, -0.55)$ on $\mathrm{BE}_2^F$ shows the boundary equilibrium $p_2$ with $W^{ss}(p_2)$ and the admissible equilibrium $p_1 \in R_1$ with $W^{ss}(p_1)$. Panel~(b) for $(\mu,\eta)=(0.0209, -0.7)$ on $\mathrm{BE}_1^F$ similarly shows the boundary equilibrium $p_1$ with $W^{ss}(p_1)$ and the admissible equilibrium $p_2 \in R_2$ with $W^{ss}(p_2)$.} 
	\label{fig:non_smooth_fold}
\end{figure}

Finally, the boundary equilibrium bifurcations along segments  $\mathrm{BE}_2^F$ and $\mathrm{BE}_1^F$ are encountered in the transition from region $\mathrm{IV}$ via region $\mathrm{VIII}$ to region $\mathrm{VI}$. In the phase portrait for $\mathrm{BE}_2^F$ in Figure~\ref{fig:non_smooth_fold}(a), the attracting boundary-node $p_2$ bounds the repelling sliding segment $\Sigma_s^r$ on the left, while a visible quadratic tangency point $F_1$ bounds it on the right. There are no pseudo-equilibria on $\Sigma_s^r$, and the equilibrium $p_1 \in R_1$ is admissible and also attracting. Trajectories above and including the union of $W^{ss}(p_2)$ in $R_2$, $\Sigma_s^r$ and the arriving orbit to $F_1$ in $R_1$ converge to $p_2$, and orbits below this union converge to $p_1$. The phase portrait along segments  $\mathrm{BE}_1^F$ in Figure~\ref{fig:non_smooth_fold}(b) is effectively the same with the roles of $p_1$ and $p_2$ exchanged. Here, the orbits below and including the union of the arriving orbit to $F_2$ in $R_2$, $\Sigma_s^r$ and $W^{ss}(p_1)$ in $R_1$ converge to the attracting boundary node $p_1 \in \Sigma_s^r$, while orbits above this union converge to the attracting equilibrium $p_2 \in R_2$.

\section{Bifurcation analysis of the smooth model}
\label{section:smooth}

We now investigate the smooth model~\eqref{eq:adjusted_welander_non_dim} for small $\varepsilon > 0$. Here we again fix the vertical mixing coefficients to $\kappa_1=0.1$ and $\kappa_2=1.0$, to enable a direct comparison of the bifurcation diagram of system~\eqref{eq:adjusted_welander_non_dim} with that of the limiting case of system~\eqref{model:PWS}. 

\begin{figure}[ht!]
	\centering
	\includegraphics{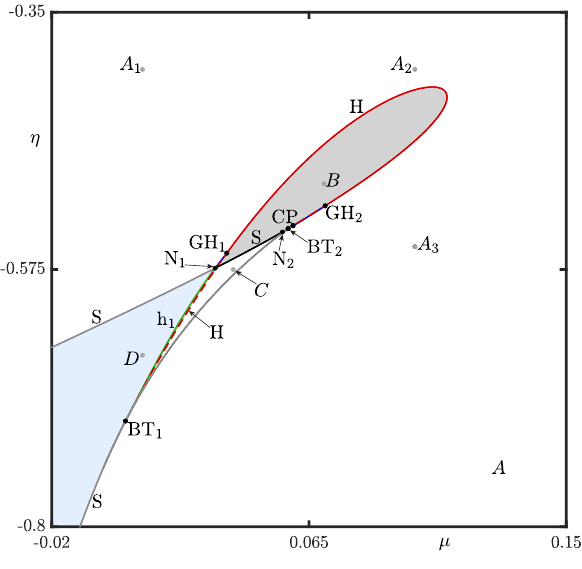}
	\caption{Two-parameter bifurcation diagram in the $(\mu,\eta)$-plane of system~\eqref{eq:adjusted_welander_non_dim} for $\varepsilon = 0.1$ and with $\kappa_1=0.1$, $\kappa_2=1.0$. Shown are curves of Hopf bifurcation $\mathrm{H}$ (red, solid when supercritical, dashed when subcritical), saddle-node bifurcation $\mathrm{S}$ (black when on periodic orbit and grey otherwise) and homoclinic bifurcation $\mathrm{h}_1$ (green), which are the main curves that divide the $(\mu,\eta)$-plane into the large regions $A, B, C$ and $D$. Also shown are codimension-two points $\mathrm{CP}$, $\mathrm{BT_1}$, $\mathrm{BT_2}$, $\mathrm{GH_1}$, $\mathrm{GH_2}$, $\mathrm{N_1}$ and $\mathrm{N_2}$; grey shading indicates the existence of a stable periodic orbit, and blue shading bistability between equilibria.}
	\label{fig:smooth_bifurcation_weak}
\end{figure}

We first consider the bifurcation diagram in the $(\mu,\eta)$-plane of system~\eqref{eq:adjusted_welander_non_dim} for the fixed value of $\varepsilon = 0.1$. It is shown in Figure~\ref{fig:smooth_bifurcation_weak} and was obtained by computing the shown bifurcation curves and codimension-two points with the continuation package AUTO-07p \cite{auto_cont}, guided by established bifurcation theory \cite{kuznetsov1998elements}. One clearly observes four main open regions, denoted $A, B, C$ and $D$, on which we focus here; associated phase portraits are shown in Figures~\ref{fig:smooth_flat_pp1}--\ref{fig:smooth_flat_pp3}.

A main element of the bifurcation diagram in Figure~\ref{fig:smooth_bifurcation_weak} is a curve $\mathrm{S}$ of saddle-node bifurcation with two branches that meet at the cusp point $\mathrm{CP}$. Along each branch of $\mathrm{S}$ there are points $\mathrm{BT}_1$ and $\mathrm{BT}_2$ of Bogdanov-Takens bifurcation (one close to $\mathrm{CP}$). From these points a curve $\mathrm{H}$ of Hopf bifurcation emerges, which is the second main element of the bifurcation diagram. Together, the curves $\mathrm{S}$ and $\mathrm{H}$ effectively form the boundaries of the four main regions $A$, $B$, $C$ and $D$. Additional ingredients are: the change of criticality of $\mathrm{H}$ at generalised Hopf points $\mathrm{GH_1}$ and $\mathrm{GH_2}$; a curve $\mathrm{h}_1$ of homoclinic bifurcation; and a segment of $\mathrm{S}$, bounded by points $\mathrm{N_1}$ and $\mathrm{N_2}$ of non-central homoclinic bifurcation \cite{strogatz_2000}, where the saddle-node bifurcation occurs on a periodic orbit (also known as SNIC or SNIPER). We remark that the complete bifurcation diagram in the $(\mu,\eta)$-plane involves
subtle additional bifurcation phenomena near the points $\mathrm{GH_1}$, $\mathrm{GH_2}$, $\mathrm{N}_1$ and $\mathrm{N}_2$ that are indistinguishable on the scale of Figure~\ref{fig:smooth_bifurcation_weak}; these include very narrow regions bounded by additional curves of homclinic bifurcation and of saddle-node bifurcation of periodic orbits, and their discussion is beyond the scope of this paper.

\subsection{Phase portraits in the main regions of the $(\mu,\eta)$-plane}
\label{section:smooth_bd_main}

\begin{figure}[t!]
\centering
\includegraphics{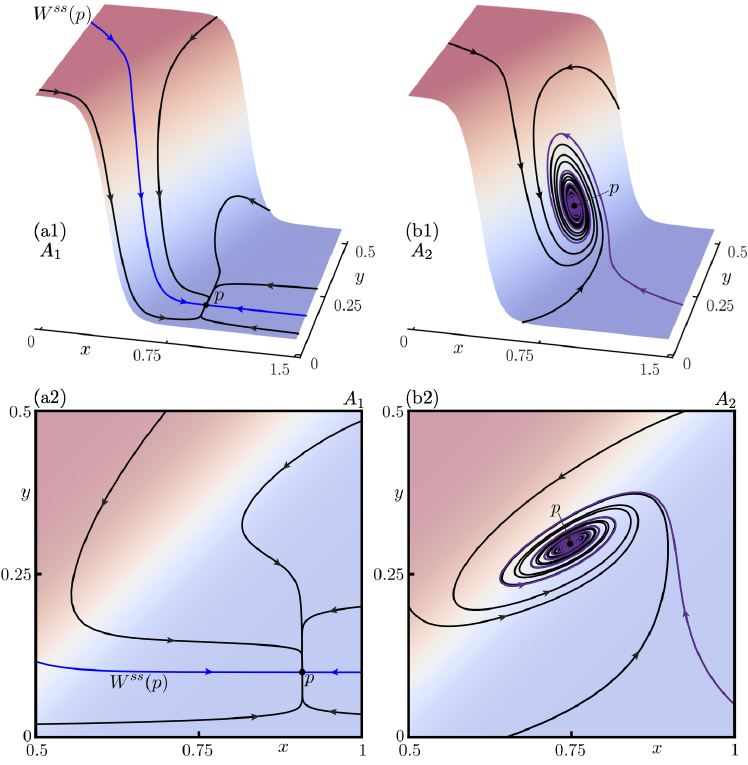}
\caption{Phase portraits at the points $A_1$ at $\mu=0.01$ and $A_2$ at $\mu=0.1$ with $\eta=-0.4$ from region $A$ in Figure~\ref{fig:smooth_bifurcation_weak}. Panels~(a1) and~(b1) shows the phase portrait on the graph of $\mc{H}_{0.1}(x,y)$, and panels~(a2) and~(b2) in the $(x,y)$-plane. Featured is the equilibrium $p$, its strong stable manifold $W^{ss}(p)$ (blue curve) when it exists, and some representative trajectories (purple curves).} 
\label{fig:smooth_flat_pp1}
\end{figure}

\begin{figure}[t!]
\includegraphics{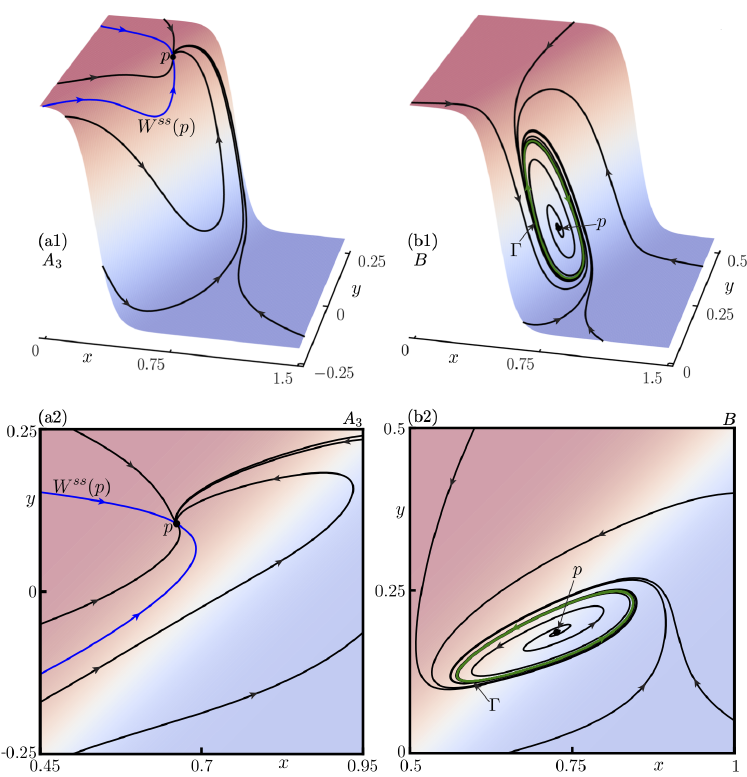}
\centering
\caption{Phase portraits at the point $A_3$ at $(\mu,\eta)=(0.1, -0.555)$ and from region $B$ at $(\mu,\eta)=(0.07,-0.50)$, shown as in Figure~\ref{fig:smooth_flat_pp1} and featuring a periodic orbit $\Gamma$ in panels~(b1) and~(b2).} 
\label{fig:smooth_flat_pp2}
\end{figure}

Region~$A$ of Figure~\ref{fig:smooth_bifurcation_weak} is bounded by the respective (supercritical) part of the curves $\mathrm{S}$ and $\mathrm{H}$. Comparison with Figure~\ref{fig:BifurcationDiag} shows that $A$ is the largest region and `covers' the five regions $\mathrm{I}$, $\mathrm{II}$, $\mathrm{III}$, $\mathrm{IV}$ and $\mathrm{VI}$ of the PWS system~\eqref{model:PWS}. Throughout region~$A$, there is a single attracting equilibrium, denoted $p$, which may correspond to distinct mixing states: weak (non-convective) mixing near $\kappa_1$, an intermediate state in between convective and non-convective mixing, or strong (convective) mixing near $\kappa_2$. This is illustrated in Figures~\ref{fig:smooth_flat_pp1} and~\ref{fig:smooth_flat_pp2}(a) with phase portraits at the parameter points labeled $A_1$, $A_2$ and $A_3$ in Figure~\ref{fig:smooth_bifurcation_weak}. We show all phase portrait of system~\eqref{eq:adjusted_welander_non_dim} in two ways to indicate when the dynamics corresponds to $\kappa_1$ or $\kappa_2$: on the graph of $\mc{H}_\varepsilon(y - x - \eta)$ over the $(x,y)$-plane and on the $(x,y)$-plane itself, where we use coloring as in Figure~\ref{fig:intro_fig2}. At parameter point $A_1$ as in Figure~\ref{fig:smooth_flat_pp1}(a), the single stable equilibrium $p$ lies in the region with $\mc{H}_{0.1}(y - x - \eta)$ near $0$ (that is, the dynamics of system~\eqref{eq:adjusted_welander_non_dim} is near $\kappa_1$), and it has real eigenvalues and a strong stable manifold $W^{ss}(p)$; hence, $p$ corresponds here to the equilibrium $p_1 \in R_1$ from regions $\mathrm{I}$ and $\mathrm{IV}$ of the PWS system~\eqref{model:PWS}. Moving to parameter point $A_2$ as in Figure~\ref{fig:smooth_flat_pp1}(b), the equilibrium $p$ now lies in the transition region where the graph of $\mc{H}_{0.1}(y - x - \eta)$ is steep; moreover, it is an attracting focus with complex conjugate eigenvalues. Finally, at parameter point $A_3$ as in Figure~\ref{fig:smooth_flat_pp2}(a), the attracting point $p$ has again real eigenvalues and a strong stable manifold $W^{ss}(p)$, and now lies in the region of the phase plane with $\mc{H}_{0.1}(y - x - \eta)$ near $1$ (that is, the dynamics is now near $\kappa_2$). Hence, $p$ now corresponds to the equilibrium $p_2 \in R_2$ in either regions $\mathrm{III}$ and $\mathrm{VI}$ of the PWS limit. We conclude that the gradual transition from $A_1$ to $A_3$ within region $A$ is very reminiscent of that from region $\mathrm{I}$, via region $\mathrm{II}$, to region $\mathrm{III}$ of system~\eqref{model:PWS}; compare with Figure~\ref{fig:pws_pp1}. 

Region~$B$ is bounded by the supercritical part of the curve $\mathrm{H}$ and the SNIPER-part of $\mathrm{S}$, and it is the `smooth version' of region $\mathrm{V}$ of system~\eqref{model:PWS}. The phase portrait in Figure~\ref{fig:smooth_flat_pp2}(b), at the marked parameter point in Figure~\ref{fig:smooth_bifurcation_weak}, shows that in region~$B$ there is indeed a stable periodic orbit $\Gamma$ surrounding the now unstable equilibrium $p$. Observe that $\Gamma$ `lives' in the switching region; that is, it lies on the steep part of the graph of $\mc{H}_{0.1}(y - x - \eta)$. Note further that the periodic orbit $\Gamma$ bifurcates at the supercritical part of the Hopf bifurcation curve $\mathrm{H}$ from the attracting focus $p$ of the phase portrait at $A_2$ in Figure~\ref{fig:smooth_flat_pp1}(b). As $\eta$ is decreased within region~$B$, the periodic orbit $\Gamma$ grows and develops two segments that lie in the region with $\mc{H}_{0.1}(y - x - \eta)$ near $0$ and near $1$, respectively; these segments correspond to the two segments of the segments periodic orbit $\Gamma$ in $\mathrm{V}$ of the limiting PWS system~\eqref{model:PWS} in Figure~\ref{fig:pws_pp2}(b).

\begin{figure}[t!] 
\includegraphics{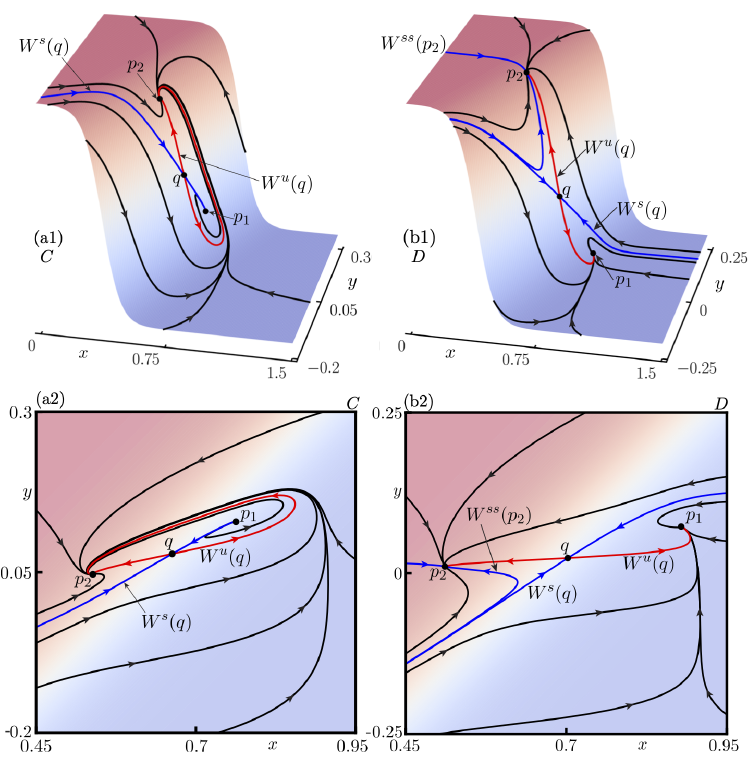}
\centering
\caption{Phase portraits in the region $C$ at $(\mu,\eta)=(0.04, -0.575)$ and region $D$ at $(\mu,\eta) = (0.01, -0.65)$. Shown in the same manner as in Figure~\ref{fig:smooth_flat_pp1}, now featuring equilibria $p_1,p_3$, and $q$ with the stable manifold $W^{s}(q)$ and unstable manifold $W^{u}(q)$.} 
\label{fig:smooth_flat_pp3}
\end{figure}

Region~$C$ of Figure~\ref{fig:smooth_bifurcation_weak} is bounded by segments of the two branches of $\mathrm{S}$ and by the homoclinic bifurcation curve $\mathrm{h}_1$ (which follows closely a subcritical part of the curve $\mathrm{H}$). Figure~\ref{fig:smooth_flat_pp3}(a) shows the representative phase portrait at the marked parameter point in Figure~\ref{fig:smooth_bifurcation_weak}. There is an attracting equilibrium, labeled $p_2$, with a high value of the transition function $\mc{H}_{0.1}$, as well as a saddle-equilibrium $q$ and a repelling equilibrium $p_1$ with an intermediate value of $\mc{H}_{0.1}$. Note that $q$ has the stable manifold $W^{s}(q)$ and unstable manifold $W^{u}(q)$, which converges to $p_2$. Region $C$ is the `smooth version' of region $\mathrm{VII}$ of the limiting PWS system~\eqref{model:PWS} in the following way: $p_2$ of the smooth system~\eqref{eq:adjusted_welander_non_dim} corresponds to $p_2 \in R_2$, and the equilibria $q$ and $p_1$ correspond to the pseudo-saddle-equilibrium $q^{-}$ and pseudo-equilibrium $q^{+}$ on the repelling sliding segment $\Sigma^r_s$, respectively; compare with Figure~\ref{fig:pws_pp3}(a).

Finally, region~$D$ is bounded by the other segments of the two branches of $\mathrm{S}$ and the homoclinic bifurcation curve $\mathrm{h}_1$. As the representative phase portrait in Figure~\ref{fig:smooth_flat_pp3}(b) at the marked parameter point in Figure~\ref{fig:smooth_bifurcation_weak} shows, it is the region of bistability and corresponds to region $\mathrm{VIII}$ of the system~\eqref{model:PWS}. The attractor $p_2$ in Figure~\ref{fig:smooth_flat_pp3}(b) is still at a high value of $\mc{H}_{0.1}$, and the saddle-equilibrium $q$ is unchanged. However, in contrast to region~$C$, the equilibrium $p_1$ is at lower value of the transition function and, moreover, it is now an attractor. The lower branch of $W^{u}(q)$ converges to $p_1$ and its upper branch to $p_2$, meaning that the stable manifold $W^{s}(q)$ forms the boundary between the basins of attraction of $p_1$ and $p_2$; compare with Figure~\ref{fig:pws_pp3}(b).

\subsection{Partial bifurcation analysis in $(\mu,\eta,\varepsilon)$-space}
\label{section:3d_par}

\begin{figure}[ht!]
	\centering
	\includegraphics{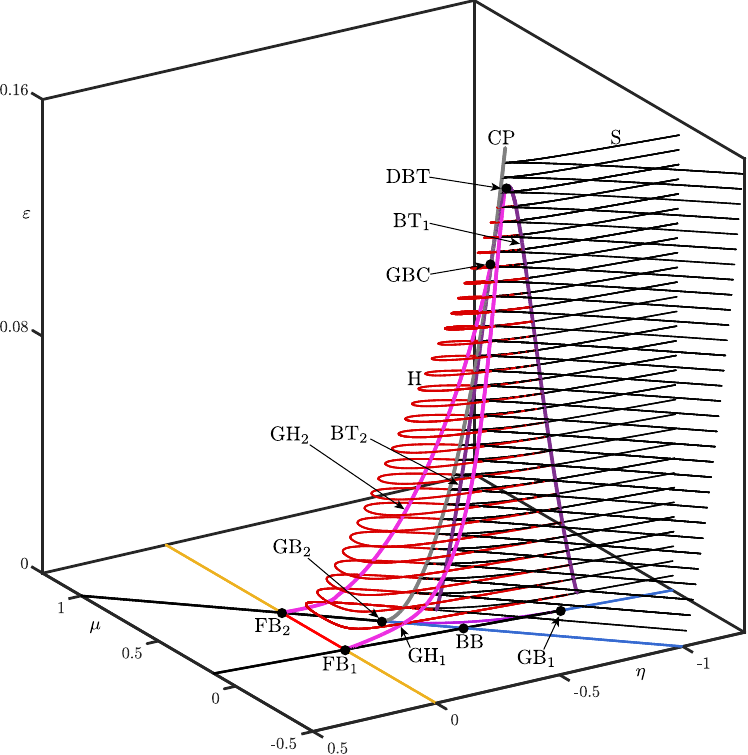}
	\caption{Partial three-parameter bifurcation diagram in $(\mu,\eta,\varepsilon)$-space of system~\eqref{eq:adjusted_welander_non_dim} for $\varepsilon>0$ and system~\eqref{model:PWS} at $\varepsilon =0$, with $\kappa_1=0.1$ and $\kappa_2=0.1$. Represented are curves $\mathrm{S}$ of saddle-node bifurcation (black) and $\mathrm{H}$ of Hopf bifurcation (red) for fixed values of $\varepsilon > 0$, together with the bifurcation diagram for $\varepsilon =0$ from Figure~\ref{fig:BifurcationDiag}. The diagram also illustrates the curve $\mathrm{CP}$ of cusp bifurcation (grey), along with branches $\mathrm{BT}_1$ and $\mathrm{BT}_2$ of Bogdanov-Takens bifurcation (dark purple) that meet at the point $\mathrm{DBT}$. Included are also curves $\mathrm{GH_1}$ and $\mathrm{GH_2}$ of generalised Hopf bifurcation (pink).} 
	\label{fig:3par}
\end{figure}

We now describe how the main elements of  the bifurcation diagram of system~\eqref{eq:adjusted_welander_non_dim} in the $(\mu,\eta)$-plane change with the switching-time parameter $\varepsilon$. Our focus here is on the curves $\mathrm{S}$ and $\mathrm{H}$, which meet at the Bogdanov-Takens points $\mathrm{BT}_1$ and $\mathrm{BT}_2$ and effectively delimit the two main regions of interest, namely region~$B$ characterised by stable oscillations and region~$D$ exhibiting bistability. Figure~\ref{fig:3par} presents the partial three-parameter bifurcation diagram in $(\mu,\eta,\varepsilon)$-space for $\varepsilon \in [0,0.16]$ and ranges of $\mu$ and $\eta$ as in Figure~\ref{fig:BifurcationDiag}. Specifically, Figure~\ref{fig:3par} shows the bifurcation diagram of the limiting system~\eqref{model:PWS} for $\varepsilon = 0$ together with the curves $\mathrm{S}$ and $\mathrm{H}$ of system~\eqref{eq:adjusted_welander_non_dim} as computed for 31 equidistant slices of fixed $\varepsilon > 0$. In this way, the corresponding surfaces $\mathrm{S}$ and $\mathrm{H}$ of saddle-node and Hopf bifurcation are visualised in $(\mu,\eta,\varepsilon)$-space with a `see-through effect'. Also shown in Figure~\ref{fig:3par} is the curve $\mathrm{CP}$ of cusp bifurcation, and the curves $\mathrm{BT}_1$ and $\mathrm{BT}_2$ of Bogdanov-Takens bifurcation, along which the surface $\mathrm{H}$ ends on the surface $\mathrm{S}$. These codimension-two bifurcation curves were computed directly by numerical continuation in $(\mu,\eta,\varepsilon)$-space. In particular, this shows that $\mathrm{BT}_1$ and $\mathrm{BT}_2$ form a single curve with a maximum at the point $\mathrm{DBT}$ at $\varepsilon \approx 0.147$, which we identified as a codimension-three degenerate Bogdanov-Takens point of focus type \cite{krauskopf2016codimension, dumortier2006bifurcations}. Additionally, curves $\mathrm{GH_1}$ and $\mathrm{GH}_2$ of generalised Hopf bifurcation are shown in Figure~\ref{fig:3par}; they were found by identifying the points the corresponding bifurcation points on the curves of Hopf bifurcation in the individual slices for fixed $\varepsilon$. We observe for increasing $\varepsilon$ that the curve $\mathrm{GH_1}$ ends at the point $\mathrm{DBT}$. The curve $\mathrm{GH_2}$, on the other hand, terminates where the curves $\mathrm{CP}$ and $\mathrm{BT_2}$ intersect at the codimension-three point $\mathrm{GBC}$ at $\varepsilon \approx 0.1208$. Similarly, the codimension-two non-central homoclinic bifurcation $\mathrm{N_1}$ and $\mathrm{N_2}$ (not shown in Figure~\ref{fig:3par}) are found to vanish as $\varepsilon$ increases, prior to $\varepsilon$ reaching the value $\varepsilon \approx 0.147$ of the point $\mathrm{DBT}$. The disappearance of $\mathrm{N_1}$ and $\mathrm{N_2}$ involves a sequence of codimension-three bifurcations, whose analysis is beyond the scope of this paper. 

We first consider the relevance of the PWS limiting system~\eqref{model:PWS} for the bifurcation diagram of the smooth system~\eqref{eq:adjusted_welander_non_dim}. While the continuation of Hopf and saddle-node bifurcation of system~\eqref{eq:adjusted_welander_non_dim} becomes very challenging for small values of $\varepsilon$ near $0$, we managed to compute the respective curves $\mathrm{S}$ and $\mathrm{H}$ in the slice at $\varepsilon=0.005$. As illustrated in Figure~\ref{fig:3par}, this turns out to be sufficient for determining the convergence of $\mathrm{S}$ and $\mathrm{H}$ to the corresponding non-smooth bifurcation as $\varepsilon$ approaches $0$. Specifically, the lower boundary of the surface $\mathrm{S}$ of saddle-node bifurcation in the $(\mu,\eta)$-plane at $\varepsilon=0$ is the union of the non-smooth fold boundary equilibrium bifurcation curve segments $\mathrm{BE}_1^F, \mathrm{BE}_2^F$ and $\widehat{\mathrm{BE}}^F_2$, and the pseudo-saddle-node bifurcation $\mathrm{PS}$. The surface $\mathrm{H}$ of Hopf bifurcation has as its  boundary at $\varepsilon = 0$ the union of the curve segments $\widehat{\mathrm{BE}}_1^P$, $\widehat{\mathrm{BE}}_2^P$ and $\widetilde{\mathrm{BE}}_1^P$ of persistence boundary equilibrium bifurcation, and $\mathrm{FU}$ of fused-focus bifurcation. Moreover, the curves $\mathrm{BT}_1$ and $\mathrm{BT}_2$ of Bogdanov-Takens bifurcation converge to the points $\mathrm{GB}_1$ and $\mathrm{GB}_2$ of generalised boundary equilibrium bifurcation, respectively; the curve $\mathrm{CP}$ of cusp bifurcation also converges to the point $\mathrm{GB}_2$. Similarly, the curves $\mathrm{GH}_1$ and $\mathrm{GH}_2$ of generalised Hopf bifurcation converge to the points $\mathrm{FB}_1$ and $\mathrm{FB}_2$ of fold-boundary equilibrium bifurcation, respectively. 

We now consider the influence of increasing the switching-time parameter $\varepsilon$. Observe from Figure~\ref{fig:3par} that the surface $\mathrm{H}$ of Hopf bifurcation `ends' at the point $\mathrm{DBT}$ at $\varepsilon \approx 0.147$. Specifically, the curve $\mathrm{H}$ in the $(\mu,\eta)$-plane for fixed $\varepsilon < 0.147$ shrinks to a point at $\mathrm{DBT}$ and disappears. Since all other curves of codimension-two bifurcations have also disappeared, above $\mathrm{DBT}$ one only finds the surface $\mathrm{S}$ of saddle-node bifurcation with the curve $\mathrm{CP}$ of cusp bifurcation. Therefore, in any slice for fixed $\varepsilon > 0.147$ the remaining regions are: region~$A$ with a single attracting equilibrium that can take any value of $\mc{H}_\varepsilon$, and the bistability region~$D$, where two stable equilibria coexist, one associated with $\mc{H}_\varepsilon$ near $0$ and the other with $\mc{H}_\varepsilon$ near $1$. In particular, both region~$B$ with stable oscillations, and region~$C$, no longer exists for $\varepsilon > 0.147$. Hence, we conclude that the existence of self-sustained oscillations in the (adjusted) Welander model requires sufficiently fast switching between convective and non-convective mixing of surface water with the deep ocean.

\section{Discussion and outlook}
\label{section:disscussion}

We studied the adjusted Welander model~\eqref{eq:adjusted_welander_non_dim} with transition function $\mc{H}_\varepsilon$ between weak and strong mixing between the warm surface and cold deep ocean as given by~\eqref{eq:transition_function}. This conceptual model in the context of the AMOC describes the evolution of temperature and salinity on the ocean surface in the Labrador and Nordic seas. We performed a bifurcation analysis with advanced tools from (non-smooth) dynamical systems theory, first for the piecewise-smooth limiting case $\varepsilon = 0$ when $\mc{H}_0$ is the Heaviside function, and then for the smooth case of $\mc{H}_\varepsilon$ with small $\varepsilon > 0$. Specifically, we presented bifurcation diagrams in the $(\mu,\eta)$-plane of salinity versus temperature flux ratio $\mu$ and density threshold $\eta$, where the rates $\kappa_1$ of weak (non-convective) and $\kappa_2$ of strong (convective) mixing were fixed at suitable values. For the PWS model with $\varepsilon = 0$, all curves of codimension-one bifurcations and points of codimension-two bifurcations were determined analytically --- resulting in a complete description of all possible dynamics and the transitions between them.  In this way, we identified the respective discontinuity-induced bifurcations, including the continuum of homoclinic orbits investigated in \cite{leifeld2016nonsmooth}, and showed how these are generated or lost as $\mu$ and $\eta$ change along different paths. In fact, the bifurcation diagram in the $(\mu,\eta)$-plane we presented for this case is complete and representative: it does not change in a qualitative way when a different choice is made for $0 < \kappa_1 < \kappa_2$, as the expressions we derived show. For the smooth case, we computed the corresponding bifurcation diagram in the $(\mu,\eta)$-plane for $\varepsilon = 0.1$ by means of numerical continuation. While the bifurcation diagram is complete, we concentrated here on four main regions of dynamics. In particular, we identified the region with oscillations found in Welander's original model \cite{welander1986thermohaline}, as well as a region of bistability that resembles previously described dynamics in a hierarchy of AMOC models \cite{weijer2019stability}. 

We also performed a partial bifurcation analysis in $(\mu,\eta,\varepsilon)$-space for small values of $\varepsilon$, which focused on surfaces of Hopf and saddle-node bifurcations that (effectively) bound the main regions. In this way, we showed how the bifurcation diagram for $\varepsilon > 0$ is `connected' to that of the PWS limit. Here, the switching time $\varepsilon$ plays the role of a parameter that desingularises the limiting Heaviside swtiching function for $\varepsilon = 0$. A direction for future mathematical work would be to use tools from geometric singular perturbation theory \cite{JonesGSPT} to study via slow-fast regularisation \cite{Jeffrey} how complicated smooth dynamics arises from the piecewise-linear limit. In this context, we conjecture that the family of homoclinic orbits along the segment $\widehat{\mathrm{BE}}_1^P$ will generate a singular Hopf bifurcation with subsequent canard explosion to the Welander-type large periodic orbit --- with the maximally sliding orbit $\gamma^*_1$ being the limit of a maximal canard \cite{slowfast_survey}.

Returning to the context of the AMOC, we found for large influxes of freshwater (smaller $\mu$) that mixing is dominantly non-convective, with the system approaching a stable equilibrium associated with $\kappa_1$. Conversely, the mixing is dominated by convection for large influxes of salinity (larger $\mu$), with convergence to a stable equilibrium associated with $\kappa_2$. We found that the intermediate region of bistability in the AMOC strength exists throughout and is rather independent of the switching time parameter $\varepsilon$.  In contrast, the region of oscillations, where the AMOC strength changes periodically between strong and weak, does depend on $\varepsilon$. In fact, oscillations are present only for sufficiently small $\varepsilon$:  when the switching between the two regimes of mixing regimes becomes too slow, oscialltions are no longer observed.

More generally, the investigation of a conceptual model, such as system~\eqref{eq:adjusted_welander_non_dim}, is a tool to uncover and highlight possible types of dynamics one may observe in the the AMOC. Specifically, we considered here the issue of deep ocean mixing in the North Atlantic in isolation from the larger climate system. Of course, there are many other climate processes that influence the overall state of the AMOC, and the analysis presented should be seen as forming a basis for the investigation of possible extensions of the model. There are several interesting directions for future research in this regard, all with their own mathematical challenges. One option is to consider additonal boxes in the model, such as an Equatorial box as in Stommel's original setup \cite{stommel1961thermohaline}, or even to model the two deep-water convection sites in the Labrador sea and the Nordic seas by separate boxed as in \cite{nkk_3box}; indeed, such models are of higher dimensions, which makes their bifurcation analysis more involved. Another direction is to incorporate seasonal changes, for example, by periodic forcing the freshwater influx parameter $\mu$, which leads to a non-autonomous model. Finally, the AMOC displays a number of feedback loops, such as the salt-advection into the subpolar North Atlantic. Incorporating feedback loops leads to the study of conceptual climate models in the form of delay differential equations, the study of which is possible but challenging because they have an infinite-dimensional phase space \cite{kkp_dde_climate}.

\section*{Acknowledgements}
This work was supported in part by Royal Society Te Ap\={a}rangi Marsden Fund grant \#19-UOA-223. We thank Henk Dijkstra for many helpful discussion, especially regarding the form of the adjusted Welander model we study here.

\appendix

\section{Proofs of Propositions~\ref{prop:escaping_sliding}--\ref{prop:pws_codim2}}

\vspace*{-2mm}
\noindent
We now state and then verify the required properties for the specific case of system~\eqref{model:PWS}, with reference to the literature on planar Filippov systems where applicable. For in-depth background on general Filippov system theory and the associated formalism see \cite{guardia2011generic,kuznetsov2003one}.

\medskip
\noindent
\textbf{Proof of Proposition~\ref{prop:escaping_sliding} (Sliding segments and tangency points).} 
	
The linear switching manifold $\Sigma$ is given as the zero set of the switching function $g(x,y) = y - x - \eta$, and it has the constant normal vector $\mathbf{n} = \binom{-1}{1}$. A tangency point $F_i$ occurs when $(f_i\cdot \mathbf{n})(x,y) = 0$ (more generally, when the first Lie derivative of $g$ with respect to $f_i$ is zero \cite{guardia2011generic}). With $y = x+ \eta$ on $\Sigma$ we obtain
	\begin{align}
\label{eq:fdotg}
		(f_i\cdot \mathbf{n})(x,x + \eta) = x - 1 + \mu - \eta \kappa_i, 
	\end{align}
which yields~\eqref{pws:tangent_points}. The visibility of the tangency point $F_i$, when it is quadratic, is determined by the curvature of the orbit of $f_i$ from $F_i$ relative to $\Sigma$. This is measured by the second Lie derivative of $g$ with respect to $f_i$ \cite{guardia2011generic}, which for system~\eqref{model:PWS} is given by 
	\begin{align*}
		(f_i \cdot \nabla (f_i \cdot \mathbf{n}))(x,y) = (1+\kappa_i)(1 - (1+\kappa_i)x) - \kappa_i\mu + \kappa_i^2y, 
	\end{align*}
	where $\nabla$ is the gradient. Evaluating at $F_i$ gives
	\begin{align*}
		(f_i \cdot \nabla (f_i \cdot \mathbf{n}))(F_i) =  \mu + \kappa_i(\mu-\eta-1) - \eta\kappa_i^2,
	\end{align*} 
which yields the genericity condition~\eqref{cond:Fi_quad}, and the visibility conditions \eqref{cond:F1_visible} and~\eqref{cond:F2_visible}.

The tangency points $F_1$ and $F_2$ bound $\Sigma_s$ and \eqref{eq:fdotg} implies
	\begin{align*}
		(f_1 \cdot \mathbf{n})(x, x + \eta) > 0\ \mathrm{for} \ x > F_1 \quad \mathrm{and} \quad (f_2 \cdot \mathbf{n})(x, x + \eta) < 0 \ \mathrm{for} \ x < F_2.
	\end{align*}
From~\eqref{pws:tangent_points} we know that $F_1 < F_2$ for $\eta > 0$, while $F_2 < F_1$ for $\eta < 0$, which yields \eqref{eq:attacking_sliding} and~\eqref{eq:repelling_sliding}.
\qed

\medskip
\noindent
\textbf{Proof of Proposition~\ref{prop:equilibria_stabilities} (Equilibria, sliding vector field and pseudo-equilibria).} 
\begin{enumerate} 
\item 
Expression~\eqref{eq:filippov_eq} immediately follow from setting $f_i(x,y) = 0$, and conditions~\eqref{cond:p1_admissible} and~\eqref{cond:p2_admissible} are immediate consequences from the definition of $R_i$ in~\eqref{eq:R1} and~\eqref{eq:R2}, respectively. The Jacobian 
		\begin{align}
\label{eq:J}
			J_{f_i}(x,y) = \left[\begin{matrix}
				-(1+\kappa_i)& 0& 
				\\
				0& -\kappa_i&
			\end{matrix}\right]
		\end{align}
of $f_i$ has two negative real eigenvalues $\lambda_{ss}= -(1+\kappa_i)$ and $\lambda_{s} =\kappa_i$, which implies that $p_i$ is a stable node.  Since $\lambda_{ss}$ has eigenvector $\binom{1}{0}$, the statement on $W^{ss}_{loc}(p_i)$ follows.
\item 
The sliding vector field on the line $\Sigma_s$ is given by
\begin{align}
\label{eq:Filipov_PWS}
	f_s(x,x+\eta) &= \left((1-\lambda(x)) f_1 + \lambda(x) f_2\right) (x,x+\eta),
\end{align}
where $\lambda(x) \in [0, 1]$ is chosen such that the vector $f_s$ is in the (constant) direction $\binom{1}{1}$ of $\Sigma_s$. This means that both components of the vector $f_s(x,x+\eta)$ are equal, which this is the case for 
\begin{align*}
	\lambda(x) = \frac{x + \mu - \kappa_2\eta - 1}{\eta(\kappa_1 - \kappa_2)}.
\end{align*}
Insertion into \eqref{eq:Filipov_PWS} and simplification yields $f_s$ as given in~\eqref{eq:extended_vector}.
\item
Setting $f_s(x,x+y) = 0$ means solving the quadratic equation 
\begin{align*}
Q(x) := \mu + (\mu - \kappa_2\eta - 1) x + x^2
\end{align*}
in~\eqref{eq:extended_vector}, which gives the expressions for $q^\pm$ in~\eqref{eq:pseudo_eq}. The stated properties follow from evaluating
$\frac{d Q(x)}{dx}$ at the $x$-values of $q^-$ and $q^+$, respectively.
\qed
\end{enumerate} 


\medskip
\noindent
\textbf{Proof of Proposition~\ref{prop:codim_1} (Codimension-one bifurcations).}
\begin{enumerate}
\item  
The equilibrium $p_i$ from Proposition~\ref{prop:equilibria_stabilities} collides with the switching manifold $\Sigma$ when
\begin{align*}
  g(p_i) = \frac{\mu}{\kappa_i} - \frac{1}{\kappa_i + 1} - \eta = 0.
\end{align*}
Solving this for $\eta$ gives the stated expression for $\mathrm{BE}_1$ and $\mathrm{BE}_2$. According to~\eqref {pws:tangent_points} the respective boundary equilibrium bifurcations happens at the tangency point $F_i = p_i$, and~\eqref{eq:pseudo_eq} shows that this involves the (dis)appearance of an admissible pseudo-equilibria through $F_i$. Simultaneously, there is a change in visibility of $F_i$ \cite{di2008discontinuity, kuznetsov2003one}, as can be seen from~\eqref{cond:F1_visible} and~\eqref{cond:F2_visible}. It follows that the visibility of the tangency points $F_i$ and the presence of admissible pseudo-equilibria $q^{\pm}$ in the different regions of the $(\mu,\eta)$-plane are as stated. See Proposition~\ref{prop:pws_codim2} for details regarding genericity conditions and different manifestations of the boundary equilibrium bifurcations along the curves $\mathrm{BE}_1$ and $\mathrm{BE}_2$.
\item 
For $\eta = 0$ we have $F^* = F_1 = F_2$ according to~\eqref{pws:tangent_points}, which is the defining property of the fold-fold bifurcation $\mathrm{FF}$; for genericity conditions and resulting different manifestations see Proposition~\ref{prop:pws_codim2}.
\item  
A saddle-node bifurcation of pseudo-equilibria occurs when the square root in~\eqref{eq:pseudo_eq} is zero, which gives 
		\begin{align*}
			\eta = -(\mu+1) + 2\sqrt{\mu}
		\end{align*}
and, hence, $\mathrm{PS}$ as stated, and also $q^{*}$ as in~\eqref{eq:qstar}. The saddle-node is generic since $\frac{d^2 Q(x)}{dx^2} = 2 \neq 0$. Since $(\mu+1)^2 - 4 \mu = (\mu-1)^2 > 0$, we know that $\eta < 0$ along the curve $\mathrm{PS}$. Hence, $q^{*}$ lies on $\Sigma_s^r$ with $F_2 < q^{*} < F_1$, and the stated bounds for $\mu$ follow.
\qed
\end{enumerate}


\medskip
\noindent
\textbf{Proof of Proposition~\ref{prop:pws_codim2} (Codimension-two bifurcations).}
	\begin{enumerate}
\item 
The expressions for $\mathrm{FB}_i$ follow immediately from Proposition~\ref{prop:codim_1} by requiring that the curve $\mathrm{FF}$ intersects the curves $\mathrm{BE}_1$ and $\mathrm{BE}_2$, respectively, yielding $p_i = F^*$. Note that these curves intersect transversely at $\mathrm{FB}_i$, and the genericity conditions for $\mathrm{FF}$ and $\mathrm{BE}_i$ are satisfied, which means that the fold-boundary equilibrium bifurcations are generic; see \cite{della2012generalized}. 

The points $\mathrm{FB}_i$ divide the fold-fold curve $\mathrm{FF}$ into segments $\mathrm{FF}_i$, where $f_1(F^{*})$ and $f_2(F^{*})$ are colinear, and a segement $\mathrm{FU}$ where they are not. With 
\begin{align*}
			 f_i \cdot \nabla (f_i \cdot \mathbf{n})(F^{*})  = \mu + (\mu-1) \kappa_i, 		
\end{align*}
we conclude that along $\mathrm{FF}_i$ the fold-fold bifurcation is for a visible and and invisible quadratic tangency, which is exactly the case VI$_1$ described in \cite{kuznetsov2003one}. It also follows that along $\mathrm{FU}$ the fold-fold bifurcation is for two invisible quadratic tangencies, and with nearby flows in opposite directions; this identifies this case as a fused-focus bifurcation according to \cite{kuznetsov2003one}. The bifurcating (crossing) periodic orbit $\Gamma$ is stable as demonstrated by the phase portraits presented in Section~\ref{section:generic_phaseportraits}. We remark that the stability of $\Gamma$ can be determined by considering the (local) return map around $F^*$ \cite{guardia2011generic, kuznetsov2003one}, but this is beyond the scope of this paper.

The point $\mathrm{FB}_i$ also divide $\mathrm{BE}_i$ locally as stated; this follows from the change of stability of the sliding segement and the associated change from $F_1 < F_2$ for $\eta > 0$ to $F_2 < F_1$ for $\eta < 0$; see Proposition~\ref{prop:escaping_sliding} and the illustrated and discussed in depth in Section~\ref{section:pp_codim1}. 
\item  
At the point $\mathrm{BB}$ of double-boundary equilibrium bifurcation there are boundary equilibrium bifurcations similtaneously at $p_1 \neq p_2$, and its location is readily found by equating expressions in Proposition~\ref{prop:codim_1} for the curves $\mathrm{BE_1}$ and $\mathrm{BE_2}$, which intersect transversally.  It follows that the division of the curves $\mathrm{BE}_i$ as are stated; this is illustrated and discussed in Section~\ref{section:pp_codim1}.
\item 
The point $\mathrm{GB}_i$ is found by equating the expressions for the curves $\mathrm{BE}_i$ and $\mathrm{PS}$ from Proposition~\ref{prop:codim_1}. Whether the boundary equilibrium bifurcation $\mathrm{BE}_i$ is of non-smooth fold or persistence type depends on the sign of the higher-order term \cite{della2012generalized}
\begin{align*}
	(\mathbf{n}\cdot(J_{f_{j}})^{-1}\cdot f_i)(p_2) = \frac{1}{\kappa_j} \left(\frac{{\kappa_i}^2}{(\kappa_i+1)^2} - \mu \right)
\end{align*}
Here $j\neq i \in \{1,2\}$ is the respective other index and $J_{f_{j}}$ is the Jacobian from \eqref{eq:J}. Hence, a sign change for the curve  $\mathrm{BE}_i$ happens at the point $\mathrm{GB}_i$; specifically, $\mathrm{BE}_i$ is of persistence type for $\mu > \frac{\kappa_i^2}{(\kappa_i+1)^2}$ and of non-smooth fold type for $\mu < \frac{\kappa_i^2}{(\kappa_i+1)^2}$. 
\qed
\end{enumerate}


\section{Phase portraits at codimension-two bifurcations} 
\label{section:pp_codim} 

\begin{figure}[ht!]
	\centering
	\includegraphics{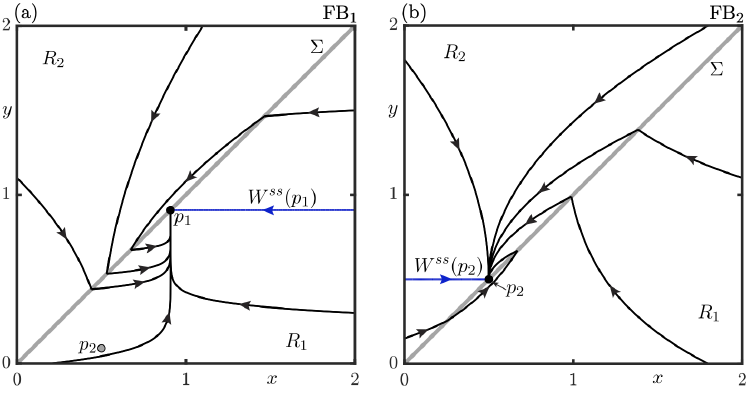}
	\caption{Representative phase portrait at codimension-two points $\mathrm{FB}_1$ and $\mathrm{FB}_2$. Panel~(a) for $(\mu, \eta) = (0.0909, 0)$ at $\mathrm{FB}_1$ shows the globally stable boundary-node $p_1$ with strong stable manifold $W^{ss}(p_1)$. Panel~(b) for $(\mu, \eta) = (0.5, 0)$ at $\mathrm{FB}_2$ shows the globally stable boundary-node $p_2$ with the strong stable manifold $W^{ss}({p_2})$. } 
	\label{fig:codim2_1}
\end{figure}

\begin{figure}[ht!]
	\centering
	\includegraphics{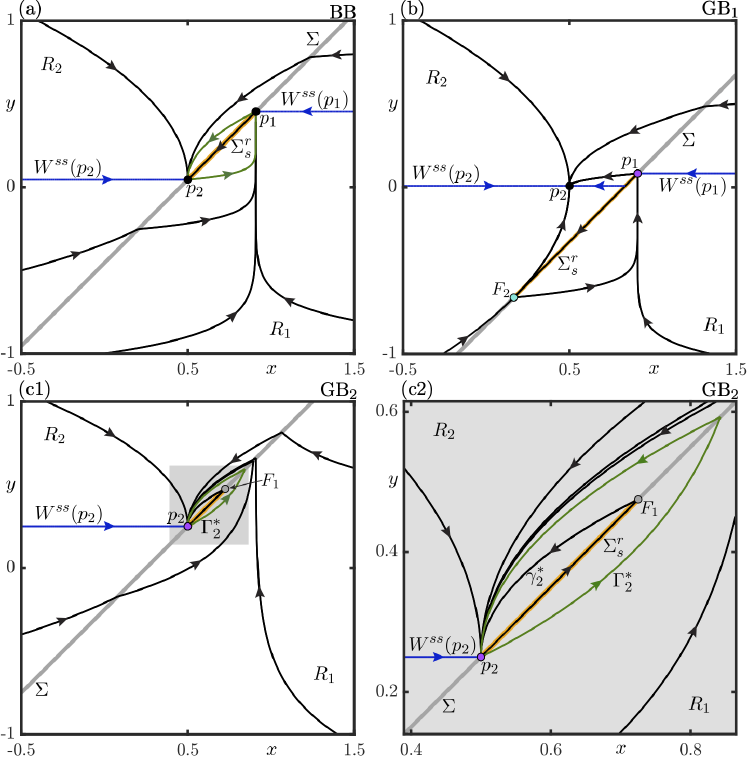}
	\caption{Representative phase portraits at codimension two points $\mathrm{GB_2}, \mathrm{BB}$ and $\mathrm{GB_1}$. Each phase portrait features a repelling sliding segment $\Sigma_s^r$. Panel~(a) for $(\mu, \eta) = (0.0083,-0.8264)$ at $\mathrm{BB}$ shows the boundary-nodes $p_1$ and $p_2$ with corresponding strong stable manifolds $W^{ss}(p_1)$ and $W^{ss}(p_2)$. A heteroclinic connection is shown in green. Panel~(b) for $(\mu,\eta) = (0.0455, -0.4545)$ at $\mathrm{GB}_1$ shows the quadratic tangency point $F_2$ and a generalized boundary-node $p_1$. The equilibrium $p_2$ is shown in $R_2$. Panel~(c1) for $(\mu, \eta) = (0.25,  -0.25)$ at $\mathrm{GB_2}$ and a magnification (c2) near the sliding segment shows an invisible quadratic tangency point $F_1$ and a generalized boundary-node $p_2$ with a strong stable manifold $W^{ss}(p_2)$. There are homoclinic connections to $p_2$ denoted $\Gamma^*_2$ and $\gamma^*_2$. } 
	\label{fig:codim2_2}
\end{figure}

We present in Figures~\ref{fig:codim2_1} and~\ref{fig:codim2_2} phase portraits at the points $\mathrm{FB_1, FB_2}$, $\mathrm{BB, GB_1}$ and $\mathrm{GB}_2$ from Proposition~\ref{prop:pws_codim2}. This illustrates how these codimension-two bifurcation points give the nearby codimension-one boundary equilibrium bifurcations $\mathrm{BE}_i$ and fold-fold bifurcations $\mathrm{FF}$ their different flavours.  

Figure~\ref{fig:codim2_1} presents phase portraits at the fold-boundary equilibrium bifurcation points $\mathrm{FB_1}$ and $\mathrm{FB_2}$. The phase portrait at $\mathrm{FB}_1$ is shown in panel~(a). It features an attracting boundary-node $p_1$ that is simultaneously a singular tangency point, which is invisible for $f_2$. All orbits converge to $p_1$ along the weak eigendirection in $R_1$. The equilibrium $p_2$ is in $R_1$ and non-admissible. The  phase  portrait at $\mathrm{FB}_2$ is shown in panel~(b). The  equilibrium  $p_2$ is now the attracting boundary-node and an invisible tangency point for $f_1$ that attracts all orbits along the weak eigendirection in $R_2$. In both cases, the strong stable manifold $W^{ss}(p_i)$ of $p_i$ is the corresponding arriving orbit in $R_i$. 

Figure~\ref{fig:codim2_2} presents phase portraits at the remaining codimension two points. The phase portrait at the double boundary equilibrium bifurcation $\mathrm{BB}$ at the intersection of the $\mathrm{BE_1}$ and $\mathrm{BE}_2$ curves is shown in panel~(a). It features a repelling sliding segment $\Sigma_s^r$ bounded on the left by the attracting boundary node $p_2$ and on the right by the attracting boundary node $p_1$. The pseudo-equilibria $q^{-}$ and $q^{+}$ are both on $\Sigma_c$ and non-admissible (and not shown): the pseudo-equilibrium $q^{-}$ is at the left hand boundary $p_2$, and $q^{+}$ is at the right hand boundary $p_1$. There is a heteroclinic connection between $p_1$ and $p_2$ composed of orbit segments in $R_1$ and $R_2$, respectively. Moreover, a (sliding) heteroclinic connection between the equilibria is composed of the sliding orbit from $p_1$ to $p_2$. If we interpret the departing orbits from $\Sigma_s^r$ as having a sliding component, then there is a continuum of homoclinic connections to $p_1$ in $R_1$ composed of departing orbits from $\Sigma_s^r$ and a corresponding sliding component. Note that the boundary-node $p_1$ has a strong stable manifold $W^{ss}(p_1)$ composed of a horizontal component in $R_1$. Similarly, the boundary-node $p_2$ also has a strong stable manifold $W^{ss}(p_2)$, composed of a horizontal component in $R_2$. Overall, the boundary equilibrium bifurcations occurring simultaneously leads to a pseudo-equilibrium emerging on the sliding segment along $\widetilde{\mathrm{BE}}_1^P$ and $\widehat{\mathrm{BE}}_2^F$; see Section~\ref{section:pp_codim1}. The phase portrait at the generalised boundary equilibrium bifurcation $\mathrm{GB}_1$ is shown in panel~(b). It features a repelling sliding segment $\Sigma_s^r$ bounded on the left by the visible quadratic tangency point $F_2$ and on the right by the attracting generalised boundary-node $p_1$ (shown in magenta). The pseudo-equilibria $q^{-}$ and $q^{+}$ undergo a pseudo-saddle-node bifurcation at the right hand boundary-node $p_1$. The phase portrait at the generalised boundary equilibrium bifurcation $\mathrm{GB}_2$ is shown in panel~(c1) with a magnification near the sliding segment in panel~(c2). The repelling sliding segment $\Sigma_s^r$ is bounded on the left by generalised boundary node $p_2$ (shown in magenta) and on the right by the invisible quadratic tangency point $F_1$. The homoclinic connections $\gamma^*_2$ and $\Gamma_2^*$ are the same as described in Section~\ref{section:pp_codim1}; see also Figure~\ref{fig:persistance_1}(b). The departing orbits from $\Sigma_s^r$ together with the respective sliding component from $p_2$, form a continuum of homoclinic connections to $p_2$.  In particular, within $\gamma^*_2$ there are homoclinic connections composed of a departing orbit from $\Sigma_s^r$ in $R_2$ and the corresponding sliding orbit. There are also homoclinic connections inbetween $\gamma_2^*$ and $\Gamma_2^*$, which feature departing orbits from $\Sigma_s^r$ in $R_1$ that cross $\Sigma$ into $R_2$. 

\bibliography{BK_AMOC_pws_references}
\bibliographystyle{elsarticle-num} 


\end{document}